\documentclass[12pt]{amsart}
\title
[Rigidity of Sierpi\'nski carpets]
{
Quasisymmetric rigidity  of  square Sierpi\'nski carpets
}
\author{Mario Bonk
and Sergei Merenkov
}

\thanks{M.B.\ was supported by NSF grants DMS-0244421, 
DMS-0456940,  DMS-0652915, DMS-1058772, and DMS-1058283.}
\thanks{S.M.\ was  supported by NSF grants 
DMS-0244421, DMS-0400636, DMS-0653439, DMS-0703617, and DMS-1001144.}

\address{Mario Bonk\\Department of Mathematics\\University of California, Los Angeles, 
Box 951555\\
Los Angeles\\ CA 90095\\USA} 
\email{mbonk@math.ucla.edu}

\address{Sergei Merenkov\\Department of Mathematics\\University of Illinois at Urbana-Champaign\\1409 W Green St\\Urbana, IL 61801\\USA}
\email{merenkov@illinois.edu}

\date{December 17, 2010}

\usepackage{amsmath, amssymb, amsthm, latexsym}
\usepackage{graphicx}
\newcommand\C{{\mathbb C}}
\newcommand\Sph{{\mathbb S}}
\newcommand\N{{\mathbb N}}

\newcommand\Hp{{\mathbb H}}
\newcommand\Z{{\mathbb Z}}
\newcommand\R{{\mathbb R}}

\newcommand\co{\colon}
\renewcommand\:{\colon}
\renewcommand\Im{\mathop{\mathrm{Im}}}
\renewcommand\Re{\mathop{\mathrm{Re}}}

\newcommand\Ga{\Gamma}
\newcommand\ga{\gamma}
\newcommand\OC{{\widehat{\mathbb C}}}
\newcommand\no{\noindent}

\newcommand\sub{\subseteq}
\newcommand\id{\operatorname{id}}
\newcommand\ra{\rightarrow}
\newcommand\QS{\operatorname{QS}}
\newcommand\diam{\operatorname{diam}}
\newcommand\Mod{\operatorname{mod}}
\newcommand\mass{\operatorname{mass}}
\newcommand\dist{\operatorname{dist}}
\newcommand\iu{{\textbf{\textit{i}}}}

\theoremstyle{plain}
        \newtheorem{theorem}{Theorem}[section]
        \newtheorem*{theorem*}{Theorem}
        \newtheorem*{conj*}{Conjecture}
        \newtheorem{lemma}[theorem]{Lemma}
        \newtheorem{proposition}[theorem]{Proposition}
        \newtheorem{corollary}[theorem]{Corollary}

\theoremstyle{definition}
        
        \newtheorem{rem}[theorem]{Remark}
         \newtheorem{rems}[theorem]{Remarks}

\theoremstyle{remark}

\numberwithin{equation}{section}

\begin{document}


\begin{abstract}
{We prove that every  quasisymmetric self-homeo\-mor\-phism  of the standard 1/3-Sierpi\'nski
carpet $S_3$ is a Euclidean  isometry. For carpets in a more general family, the  standard 
$1/p$-Sierpi\'nski carpets $S_p$, $p\ge 3$ odd,  we show that the groups of quasisymmetric self-maps are finite dihedral. We also establish  that $S_p$ and $S_q$ are 
quasisymmetrically equivalent only if $p=q$. The main tool in the proof for these facts is a new invariant---a certain discrete modulus of a path family---that is preserved under  quasisymmetric maps of carpets.}
\end{abstract}

\maketitle

\section{Introduction}\label{S:Intro}
\begin{figure}
[htbp]
\begin{center}
\includegraphics[height=40mm]{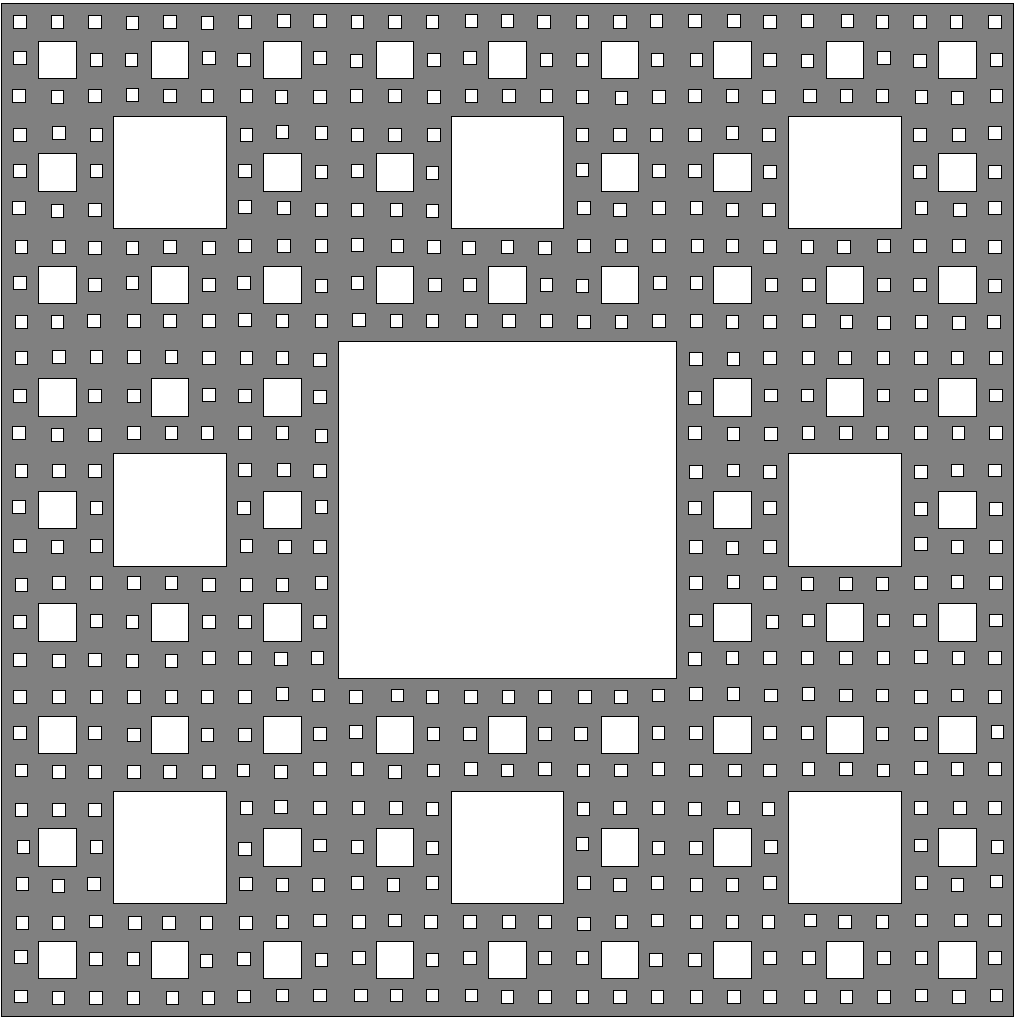}
\caption{
The standard Sierpi\'nski carpet $S_3$.
}
\label{F:Sierpinski}
\end{center}
\end{figure}

\no
The {standard Sierpi\'nski carpet} $S_3$ is a subset of
the plane $\R^2$  defined as follows. Let $$Q_0=\{(x,y)\in\R^2\co 0\leq x\leq 1,\, 0\leq y\leq1\}$$ denote the closed unit square in $\R^2$. We subdivide $Q_0$ into
$3\times3$ subsquares of equal size in the obvious way and remove the
interior of the middle square. The resulting set $Q_1$ consists of eight squares of sidelength $1/3$. 
Inductively, $Q_{n+1},\ n\geq1$, is obtained from $Q_n$ by subdividing each  of the remaining squares in the subdivision of $Q_{n}$ 
 into $3\times3$ subsquares and removing the interiors of the middle squares.  
The \emph{standard Sierpi\'nski carpet} $S_3$ is the intersection of all the sets $Q_n,\ n\geq0$ (see Figure~\ref{F:Sierpinski}).
For arbitrary  $p\geq3$  odd, the \emph{standard $1/p$-Sierpi\'nski carpet}
 $S_p$ is the subset of the plane obtained in a similar way by subdividing the  square $Q_0$ into
$p\times p$ subsquares of equal size, removing the
interior of the middle square, and repeating these operations as above. 

In general, a \emph{(Sierpi\'nski) carpet} is a metrizable topological space $S$ 
homeomorphic to the standard Sierpi\'nski carpet $S_3$. According to the
topological characterization of Whyburn~\cite{gW58}, $S$ is a carpet if and
only if it is a planar continuum of topological dimension $1$
that is locally connected and has no local cut points. 

Let  
$$\Sph^2=\{(x,y,z)\in \R^3\co x^2+y^2+z^2=1\} $$ 
denote the unit sphere in $\R^3$.  
In the following we often identify $\Sph^2$ with the extended complex plane $\OC=\C\cup\{\infty\}$ by  stereographic projection.
For subsets of $\Sph^2$ a more explicit characterization of carpets can be given as follows \cite{gW58}.  A set $S\subseteq \Sph^2$ is a carpet if and only if it can be written as 
\begin{equation}\label{eq:carprep}
S=\Sph^2\setminus\bigcup_{i\in \N} D_i,
\end{equation}
where   for each $i\in \N$ the set $D_i\subseteq \Sph^2$
is a Jordan region and the following conditions are satisfied:  $S$ has empty interior, ${\rm{diam}}(D_i)\to0$ as $i\to\infty$, and $\partial D_i\cap\partial D_j=\emptyset$ for $i\neq j$. 
This characterization implies that all the sets  $S_p$, $p\ge 3$ odd, are indeed carpets.

A  Jordan curve 
in a carpet $S$ is called a \emph{peripheral circle} if its complement in $S$ is a connected set. If $S\subseteq \Sph^2$ is a carpet, written as in \eqref{eq:carprep},   then the  peripheral circles 
 of $S$ are precisely  the boundaries $\partial D_i$ of the 
 Jordan regions $D_i$, $i\in \N$.

Let $f\co X\to
Y$ be a homeomorphism between two metric spaces $(X, d_X)$ and
$(Y, d_Y)$. The map $f$ is called \emph{quasisymmetric} if there
exists a homeomorphism $\eta\co [0,\infty)\to[0,\infty)$ such that
$$
\frac{d_Y(f(u),f(v))}{d_Y(f(u),f(w))}\leq
\eta\bigg(\frac{d_X(u,v)}{d_X(u,w)}\bigg)
$$
whenever $u,v,w\in X$, $u\ne w$.  If we want to emphasize the distortion function  $\eta$, we say that $f$ is $\eta$-\emph{quasisymmetric}.
When we speak of a quasisymmetric map $f$ from  $X$ to $Y$, then it is understood that $f$ is a homeomorphism of $X$ onto $Y$
and that the underlying metrics on the spaces have been specified.
 Unless otherwise indicated,  a carpet as in \eqref{eq:carprep}
 is equipped with the spherical metric. The carpets $S_p$ will carry the Euclidean metric.  Note that for a compact subset $K$ of $\C\sub \C\cup\{\infty\}\cong \Sph^2$ the Euclidean and the spherical metrics are comparable. So for 
 the notion of a quasisymmetric map on $K$  is does not matter which of these two metrics we choose on $K$.   
 
It is immediate that restrictions, inverses, and compositions of  quasisymmetric maps are quasisymmetric. If there is a quasisymmetric map  between two metric spaces $X$ and 
$Y$, we say that $X$ and $Y$ are  \emph{quasisymmetrically equivalent}.   The quasisymmetric self-maps on a metric space $X$, i.e., the quasisymmetric homeomorphisms of $X$ onto itself, form a group that we 
denote by $\QS(X)$. If two metric spaces $X$ and $Y$ are quasisymmetrically equivalent, then $\QS(X)$ and $\QS(Y)$ are isomorphic 
groups.

From the topological point of view all carpets are the same and so the topological universe of all  carpets consists of a single point. A much richer structure emerges if we look at metric carpets from the point of view of quasiconformal geometry. In this case, we identify two metric carpets if and only  if they are quasisymmetrically equivalent. Even if we restrict ourselves 
to carpets contained in $\Sph^2$, then the set  of 
 quasisymmetric equivalence classes of carpets is uncountable. 
  
One way to see this is to invoke a rigidity  result that has  recently been established in \cite{BKM06}. To formulate it, we call a carpet $S\sub \Sph^2$ {\em round} if its peripheral circles are geometric circles. So if $S$ is written 
as in \eqref{eq:carprep}, then each  Jordan region $D_i$ is an 
open spherical disk.  According to \cite{BKM06} two round carpets $S$ and $S'$ of measure zero are quasisymmetrically equivalent only if they are {\em M\"obius equivalent}, i.e., one is the image of the other under a M\"obius transformation on $\OC\cong \Sph^2$.  Since the group of 
M\"obius transformations depends on 6 real parameters, but 
the set of round carpets is  a family depending on essentially a countably infinite set of 
real parameters (specifying the radii and the locations of the centers  of the  complementary disks of the round carpet), it easily follows that the set of 
quasisymmetric equivalence classes of round carpets in $\Sph^2$ has the cardinality of the continuum.

Among the round carpets there is a particular class of carpets that are 
distinguished by their symmetry.  Namely, suppose that $K$ is a convex subset of hyperbolic $3$-space $\Hp^3$ with non-empty interior and non-empty  totally geodesic boundary, and suppose that there exists a group $G$  of isometries of $\Hp^3$ that leave $K$ invariant, and that $G$ acts cocompactly and properly discontinuously on $K$.
If we   identify $\Sph^2$ with the boundary at infinity $\partial_\infty \Hp^3$ of $\Hp^3$, then the limit set  $\Lambda_\infty(G)\sub \partial_\infty \Hp^3=\Sph^2$ of $G$  is a round carpet. The group 
$G$ induces an action on $\Sph^2$ by M\"obius transformations that 
leave $S=\Lambda_\infty(G)$ invariant.  Moreover, this  action is cocompact on triples of $S$.
We can  consider $G$ as a subgroup of $\QS(S)$. An immediate consequence is that $\QS(S)$ is infinite,  and that there are only finitely 
many distinct orbits of peripheral circles under the action of $G$, and 
hence of $\QS(S)$,  on $S$. In this sense, $S$ is very symmetric. 

According to an open    conjecture by Kapovich and Kleiner~\cite{KK00}, 
up to  virtual isomorphism the groups $G$ as above are precisely the  Gromov hyperbolic groups whose boundaries  at infinity are Sierpi\'nski carpets.     In order to get a better understanding  of the relevant issues in this problem, 
it seems desirable to characterize 
the carpets $S$ that arise from such groups $G$ from the point of view of their quasiconformal geometry.

To formulate these questions more precisely, let $\mathcal{S}$, $\mathcal{R}$, $\mathcal{G}$, respectively,  denote the set of all quasisymmetric equivalence 
classes of all carpets in $\Sph^2$, all round carpets, and all {\em round group carpets}, i.e., all carpets arising as limit sets $\Lambda_\infty(G)$ of groups $G$ as 
above. Then $\mathcal{G}\sub \mathcal{R}\sub \mathcal{S}$. Let  $[S]$
denote the quasisymmetric equivalence class of a carpet $S\sub \Sph^2$.  

 An obvious problem   is  where $[S_p]$ is placed in the universe $\mathcal{S}$.   It follows from the main result in \cite{mB04} that each standard carpet $S_p$ 
 is quasisymmetrically equivalent to a round carpet. Hence $[S_p]\in 
 \mathcal{R}$.  The question whether actually $[S_p]\in 
 \mathcal{G}$  arose  in discussions with B.~Kleiner and 
 the first author about ten years ago.  At the time this problem was  considered as completely inaccessible, and one stood helpless in front 
 of these and other  problems of quasiconformal geometry 
 (another well-known hard problem related to carpets is the question of the (Ahlfors regular) conformal dimension of $S_3$; see \cite{MT} for general background and  \cite{jK10} for specific results on $S_3$).


 The  main results of this paper give answers  to some of these questions.  We consider them as an important step in a better understanding of the quasiconformal geometry of 
 Sierpi\'nski carpets and hope that some of our results and techniques may be useful for progress on  the Kapovich-Kleiner conjecture. 
 
\begin{theorem}\label{T:Onethird}
Every quasisymmetric self-map of the standard Sier\-pi\'n\-ski carpet $S_3$ is a Euclidean  isometry.\end{theorem}

The isometries of $S_3$ are given  by the Euclidean symmetries 
that leave $S_3$, and also the unit square $Q_0$,  invariant. 
They form a dihedral group with $8$ elements. 
We conjecture that  also for $p>3$ each   quasisymmetric self-map of $S_p$ 
is an  isometry.  With some effort our  proof of this fact for $p=3$ can actually be extended to $p=5$,    but the general case remains open (see Remark~\ref{R:Onefifth} for more discussion).

We are able to  prove that ${\rm QS}(S_p)$ is a finite dihedral group for each odd $p$.

\begin{theorem}\label{T:Onepth}
For every odd integer $p\geq3$, the group of quasisymmetric self-maps ${\rm QS}(S_p)$ of the standard Sierpi\'nski carpet $S_p$ is  finite dihedral. 
\end{theorem}

Theorems~\ref{T:Onethird} and \ref{T:Onepth} are quite unexpected as  the group of {\em all}  homeomorphisms on $S_p$ is  large. For example, if $u$ and $v$ are two points in $S_p$ that do not 
lie on a peripheral circle of $S_p$, then there exists a homeomorphism $f\: S_p\ra S_p$ with $f(u)=v$. 

Every  bi-Lipschitz homeomorphism between   metric spaces, i.e., every  homeomorphism that distorts distances by an at most bounded multiplicative amount,  is a quasisymmetry. So Theorems~\ref{T:Onethird} and \ref{T:Onepth} remain true if one only considers 
bi-Lipschitz homeomorphisms instead of quasisymmetries. 
In general, these maps form  a rather restricted subclass of all quasisymmetries. In view of this,  one may wonder whether the bi-Lipschitz versions of Theorems~\ref{T:Onethird} and \ref{T:Onepth} are easier to establish.  Our methods do not offer any simplifications 
for this more restricted class, and it seems that there is no straightforward way to benefit from the  stronger bi-Lipschitz hypothesis on the maps.

An immediate consequence of Theorem~\ref{T:Onepth}  is that 
$[S_p]\not\in \mathcal{G}$. Indeed, $\QS(S_p)$ is a finite group, while 
$\QS(S)$ is infinite if $[S]\in \mathcal{G}$. 
So the points $[S_3],\, [S_5],\, [S_7], \dots$ lie in $\mathcal{R}\setminus \mathcal{G}$. 
As the following theorem shows, these  points are actually all distinct. 

\begin{theorem}\label{T:Standard}
Two standard Sierpi\'nski carpets $S_p$ and $S_q$, $p,q\ge 3$ odd,  are quasisymmetrically equivalent if and only if  $p=q$.
\end{theorem}
It  was previously known that if $|p-q|$ is large,  then   $S_p$ and $S_q$ cannot be quasisymmetrically equivalent; more precisely, if  $p>q$ say, and 
$$1+\frac{\log(p-1)}{\log p}>\frac{\log(q^2-1)}{\log q}, 
$$ 
then $[S_p]\ne [S_q]$. 
Here the quantity on the right  of the inequality is the Hausdorff dimension of $S_q$, while the quantity on the left  is a lower bound for the (Ahlfors regular) conformal dimension of $S_p$, i.e., for the infimum of  the Hausdorff dimensions of all Ahlfors regular metric spaces 
quasisymmetrically equivalent to $S_p$. So the inequality  guarantees that   $[S_p] \ne[S_q]$. The bi-Lipschitz version of Theorem~\ref{T:Standard} is easy to establish. Namely, if $p\ne q$, then $S_p$ and $S_q$ have different Hausdorff dimensions.  So there  cannot be any bi-Lipschitz homeomophism between these spaces, because    bi-Lipschitz maps preserve Hausdorff dimension. 

One of the main difficulties in the proof of Theorem~\ref{T:Onethird} is that we have no a priori normalization of a   quasisymmetric self-map $f$ of $S_3$. If we knew in advance, for example, that $f$ sends each  corner of the unit square to another such corner, then this  statement would immediately follow from the following theorem which is relatively easy  to establish.

\begin{theorem}\label{P:Cor}
Let $S$ and $\widetilde S$ be  square carpets of measure zero in rectangles
$K=[0,a]\times[0,1]\sub \R^2$ and $\widetilde K=[0,\tilde a]\times[0,1]\sub\R^2$,
respectively, where $a,\tilde a>0$. If $f$ is an orientation-preserving quasisymmetric homeomorphism from
$S$ onto $\widetilde S$ that takes the corners  of $K$ to the corners of
$\widetilde K$  such that $f(0)=0$, then $a=\tilde a$, $S=\widetilde S$, and
$f$ is the identity on $S$.
\end{theorem}

Here the expression {\em square carpet in a rectangle}  is used in the specific sense of the more general concept  
 of a {\em square carpet in a closed Jordan domain} 
 defined in Section~\ref{S:Corners}. A quasisymmetric map between   carpets in $\Sph^2$ is called {\em orientation-preserving} if it has 
an extension to a homeomorphism on $\Sph^2$ with this property.

Theorem~\ref{P:Cor} is analogous to the uniqueness part of~\cite[Theorem~1.3]{oS93}. Our proof is similar in spirit, but we use the classical conformal modulus instead of a discrete version of it.

Another situation where  a natural normalization implies a strong  rigidity statement is  for 
square carpets in $\C^*$-cylinders.  

\begin{theorem}\label{P:ACyl} 
Let $S$ and $\widetilde S$ be square carpets  of measure zero in  $\C^*$-cylinders $A$ and $\widetilde A$, respectively. 
Suppose that $f$ is an orientation-preserving  quasisymmetric homeomorphism of  $S$ onto $\widetilde S$ that maps the inner and outer boundary components of $A$ onto the inner and outer boundary components of $\widetilde A$, respectively. Then $f$ is (the restriction of) a  map of the form $z\mapsto f(z)=az$, where $a\in \C\setminus\{0\}$. 
\end{theorem}
See Section~\ref{S:Aux} for the relevant definitions. 
Similar rigidity results for so-called slit carpets were established 
 in~\cite{sM10}.

We now discuss some of the ideas in the proof of our main results and give a general outline of the paper. Most of the results have been announced in 
\cite{Bo2}. 

The main new tool used in proving  Theorems~\ref{T:Onethird}--\ref{T:Standard} is {\em carpet modulus}, a  version of Schramm's transboundary modulus~\cite{oS95} for path families adapted to Sierpi\'nski carpets. This is discussed in Section~\ref{S:Mod}. We also need a notion of  carpet modulus that  takes a  group action into account, see Section~\ref{S:Cmodgp}. 
The crucial feature of carpet modulus is that it is invariant under quasisymmetric maps in a suitable sense (see Lemma~\ref{L:Invmod}). 
 
 We denote by $O$  the boundary of the unit square $Q_0$ 
and by $M$ (for fixed   $p\ge 3$ odd) the boundary of the first square removed from 
$Q_0$ in the construction of $S_p$ (``the middle square"). Then 
the pair $\{O,M\}$ is distinguished by an extremality property for carpet modulus among all pairs of peripheral circles of $S_p$ (Lemma~\ref{L:Pair}). It follows that every quasiymmetric self-map $f$ of $S_p$ must preserve  the pair $O$ and $M$, i.e., $\{f(O), f(M)\}=\{O,M\}$ (Corollary~\ref{C:Group}). In principle, $f$ may interchange $O$ and $M$, but  by a  more refined analysis we will  later establish  that $f(O)=O$ and $f(M)=M$ (Lemma~\ref{L:f(O)=O}). 
This is quite in contrast to the behavior of general homeomorphisms on a carpet:  if  we have two finite  families each consisting of the same number of distinct peripheral circles of a carpet $S$, then we can find a self-homeomorphism of  $S$ that sends one family to the other family.

The proof of Corollary~\ref{C:Group} relies on some previous work.  In Section~\ref{S:Aux} we collect certain   uniformization and rigidity results that were  established  in~\cite{mB04} and~\cite{BKM06},  and derive some consequences.  Among these results is Proposition~\ref{P:Expldescr} which gives an 
  explicit description of extremal mass distributions for carpet modulus of certain path  families. This is an important ingredient in the proof of  Corollary~\ref{C:Group}.
   Corollaries~\ref{C:TCT}, \ref{C:Cycl}, \ref{C:CC}, and~\ref{C:Ccp} in Section~\ref{S:Aux}  give information on quasisymmetric maps on certain carpets under various normalizing conditions for points and peripheral circles. 
   
   In Section~\ref{S:Corners} we prove Theorems~\ref{P:Cor}
   and~\ref{P:ACyl}. This is essentially independent of the rest of the paper, but Theorem~\ref{P:Cor} will later be  used in the proof of Theorem~\ref{T:Standard}. 
  
  The fact that every   quasisymmetric self-map of $S_p$
  preserves the pair $\{O,M\}$ 
  already has some strong consequences. For example, combined with the results in Section~\ref{S:Aux}, one can easily  derive that the  group ${\rm QS}(S_p)$ is finite (Corollary~\ref{C:Gpfinite}).    
  To push the analysis further and to arrive at proofs of Theorems~\ref{T:Onethird}--\ref{T:Standard}, we need one additional essential 
  idea;  namely, we will investigate weak tangent spaces 
  of the carpets $S_p$ and induced quasisymmetric maps on these weak tangents (see Section~\ref{S:Wtangents}). In particular,  we prove that the  weak tangent of $S_p$ at a  corner of $O$ cannot be mapped  to the weak tangent  of $S_p$ at a corner of $M$ by a (suitably normalized) quasisymmetric map 
  (Proposition~\ref{P:Corwtangents}). Actually, we conjecture that such maps only exist if the weak tangents are isometric, but Proposition~\ref{P:Corwtangents} is the only result in this direction that we are able to prove. 
  
Using these statements on weak tangents, we will give proofs 
of Theorems~\ref{T:Onethird}--\ref{T:Standard}  in the following Section~\ref{S:PT1}.  Overall, the ideas in these  proofs  are very similar. In order to establish  Theorem~\ref{T:Standard}, for example,  one wants to apply  Theorem~\ref{P:Cor}. For this one essentially has to show that a quasisymmetric map  
$f\: S_p \ra S_q$ preserves the set of corners of $O$.
Let $M_p$ and $M_q$ denote the boundary of the middle square for $S_p$ and $S_q$, respectively. Using 
the  extremality property for the pair $\{O,M\}$, one  can show that 
$\{f(O), f(M_p)\}=\{ O,M_q\}$. This leads to various combinatorial 
possibilities, and in each case one analyzes what happens to the corners of $O$ under the map $f$. The case $f(O)=O$ leads to a favorable  situation, where the set of corners of $O$ is  preserved and where one can apply Theorem~\ref{P:Cor} to conclude $S_p=S_q$. One wants to rule out the existence of the map $f$ in the other cases,  for example when $f(O)=M_q$. In all these cases, one eventually ends up with a contradiction to 
Proposition~\ref{P:Corwtangents}.

\medskip\noindent 
\textbf{Acknowledgments.} The authors are indebted to Bruce Kleiner and the late Juha Heinonen for many fruitful discussions.

This work was completed while the second author was visiting 
the Hausdorff Research Institute for Mathematics, Bonn, Germany,  in the fall of 2009, and the Institute for Mathematical Sciences, Stony Brook, New York, in the fall of 2010. He thanks these institutions for their hospitality. 

\section{Carpet modulus}\label{S:Mod}
 
\no 
We first make some remarks about notation and terminology used in the rest of the paper. We denote the imaginary unit in $\C$ by $\iu$. 
Let  $(X,d)$ be  a metric space. If $x\in X$ and $r>0$,  we denote by $B(x,r)$ the open ball, and by 
$\overline B(x,r)$ the closed ball  in $X$ that has radius $r>0$ 
and is centered at $x$. 
If $\lambda>0$ and $B=B(x,r)$, we let   $\lambda B$ be the open ball of radius  $\lambda r$ centered at $x$. 
If $A\sub X$, then $\diam(A)$  is the diameter,  $\chi_A$ the characteristic function, and $\#A\in \N_0\cup\{\infty\}$ the cardinality of $A$.  If $B\sub X$ is another set, then we let  
$${\rm dist}(A,B):=\inf\{d(a,b)\co a\in A,\,  b\in B\}$$ 
be  the distance between $A$ and $B$.

If $X$ is a set, then ${\rm id}_X$ is the identity map on $X$. 
If $f\: X\ra Y$ is a map between two sets $X$ and $Y$, and $A\sub X$, then $f|A$ denotes the restriction of $f$ to $A$. 

Unless otherwise indicated,  our ambient metric space is the sphere $\Sph^2$ 
equipped with the spherical metric induced by the  standard Riemannian structure on   $\Sph^2$. In this metric space the balls
are spherical disks.

A {\em Jordan region} in $\Sph^2$ is an open connected set bounded by a Jordan curve, i.e., a set homeomorphic to a circle. A {\em closed Jordan region} in $\Sph^2$ is the closure of a Jordan region.

 A \emph{path} $\gamma$ in a metric space $X$ is a continuous map $\ga\: I\ra X$  of a finite  interval $I$,  i.e., a set of the form $[a,b],[a,b),(a,b]$, or $(a,b)$, where $a< b$ are real numbers,
into the space $X$. If $\ga$ is a map from $(a,b)$ we say that the path  is \emph{open}. As is standard, we often denote by $\ga$ also the image set $\ga(I)$ in $X$. The limits $\lim_{t\to a}\gamma(t)$ and $\lim_{t\to b}\gamma(t)$, if they exist,  are called \emph{end points} of $\gamma$.  If $A,B\sub X$ then we say that
$\ga$ {\em connects} $A$ and $B$ if $\ga$ has endpoints and one of them lies in $A$ and the other in $B$. A path is called a {\em subpath} of $\ga$ it is of the form $\ga|J$ for some interval $J\sub I$. 
We  denote the length of $\ga$ by ${\rm length}(\ga)$.  The path $\ga$ is called  {\em rectifiable} if it has finite length, and {\em locally rectifiable} if $\ga|J$ is a rectifiable path for every compact subinterval $J\sub I$. 

Let $\sigma$ denote the spherical measure and $ds$ the spherical line element 
on $\Sph^2$ induced by the standard Riemannian metric. A \emph{density} $\rho$ is a non-negative Borel  function defined on $\Sph^2$. The density $\rho$ provides a pseudo-metric with line element $\rho\, ds$. If $\Ga$ is  a family of paths in $\Sph^2$, then  the \emph{conformal modulus} of $\Ga$, denoted ${\rm{mod }}(\Ga)$, is defined to be the infimum of the \emph{mass} 
$$
\int\rho^2\, d\sigma
$$
over all \emph{admissible} densities $\rho$, i.e., all densities such that for $\rho$-\emph{length} of each  locally rectifiable path $\gamma\in \Ga$  we have the inequality
$$
\int_{\gamma}\rho\, ds\ge 1. $$
If $\rho$ is admissible and has minimal mass among all densities admissible for $\Gamma$, then $\rho$ is called {\em extremal}.   
Often it is convenient to change the spherical metric that was the underlying {\em base metric} in the definition of 
${\rm{mod }}(\Ga)$ to another conformally equivalent metric. This leads to the same quantity ${\rm{mod }}(\Ga)$ (see \cite[Remark~6.1]{mB04} for more discussion).

Conformal modulus is  
monotone~\cite[Section~4.2, p.~133]{LV}:
if  $\Gamma$ and $\Ga'$ are two path families in $\Sph^2$ such that every path $\ga\in\Ga$ contains a subpath $\ga'\in\Ga'$, then 
$$
{\rm mod}(\Ga)\leq{\rm mod}(\Ga').
$$  
In particular, this inequality holds if $\Gamma\sub \Gamma'$.

Conformal modulus is also countably subadditive~
\cite[Section~4.2, p.~133]{LV}: 
 for any countable union $\Ga=\bigcup_i\Ga_i$ of path families $\Gamma_i$ in $\Sph^2$ we have 
$$
{\rm mod}(\Ga)\leq\sum_i{\rm mod}(\Ga_i).$$
Here and in the following we adopt the convention that if the range of an index 
such as $i$ above is not specified, then it is extended over $\N$, i.e., it runs through $1,2,\dots$.

An important property of conformal modulus is its invariance under conformal and its   quasi-invariance  under quasiconformal maps. The latter means that if $\Gamma$ is a family of paths contained in a region $D\sub \Sph^2$ and   if $f\co D\to\widetilde D$ is a quasiconformal map onto another region $\widetilde D\sub \Sph^2$, then 
\begin{equation}\label{E:Qcqinv}
\frac1K{\rm{mod }}(\Ga)\leq {\rm{mod }}(f(\Ga))\leq K{\rm{mod }}(\Ga),
\end{equation}
where  $f(\Gamma):=\{f\circ \gamma:\gamma\in \Gamma\}$ and $K$ depends only on the dilatation  of $f$~\cite[Theorem~3.2, p.~171]{LV}. 
 For the basic definitions and general background on quasiconformal maps see \cite{lA66, LV, Va}. 
 We adopt the ``metric definition" of  quasiconformal maps and allow them to be orientation-reversing. 
 
 If a certain property for paths in a family $\Gamma$  holds for all paths outside an exceptional 
 family $\Ga_0\sub \Gamma$ with ${\rm mod}(\Ga_0)=0$, we say that it holds for \emph{almost every path} in $\Gamma$.

Now let $S\sub \Sph^2$ be a carpet  as in \eqref{eq:carprep}, and
$\Ga$ be  a family of paths in
$\Sph^2$. Then we  define  the \emph{carpet modulus}  of $\Ga$
(\emph{with respect to} $S$), denoted by ${\rm mod}_S(\Ga)$,  as follows. Let  $\rho$ be a \emph{mass distribution} defined on the peripheral circles of $S$, i.e., a function $\rho$ that assigns to each peripheral circle $C_i=\partial D_i$ of $S$ a non-negative number $\rho(C_i)$. 
If $\ga$ is a
path in $\Sph^2$, the $\rho$\emph{-length} of $\ga$ is
$$
\sum_{\ga\cap C_i\neq\emptyset}\rho(C_i).
$$ 
We say that a mass distribution $\rho$ is \emph{admissible} for
${\rm mod}_S(\Ga)$ if  for almost every path
$\ga\in\Ga$ the $\rho${-length} of $\ga$ is  $\ge 1$; so we require  that there exists a family $\Ga_0\sub \Gamma$ with $\Mod(\Gamma_0)=0$ such that   
$$\sum_{\ga\cap C_i\neq\emptyset}\rho(C_i)\ge 1$$ 
for  every path
$\ga\in\Ga\setminus \Gamma_0$. We call $\Gamma_0$ an  {\em exceptional  family} for $\rho$. Now we set 
$${\rm mod}_S(\Ga)=\inf_\rho\bigg\{\sum_i\rho(C_i)^2\bigg\},
$$
where the infimum is taken over all mass distributions $\rho$ that 
are admissible for ${\rm mod}_S(\Ga)$. 
The sum $\sum_i\rho(C_i)^2$ is called the \emph{(total) mass} of
$\rho$, denoted ${\rm mass}(\rho)$.  Often we will consider 
a mass distribution $\rho$ of finite mass as an element in  the Banach space 
$\ell^2$ of square summable sequences. By definition $\ell^2$ consists of all sequences $a=(a_i)$ with $a_i\in \R$ for $i\in \N$ and 
$$ \Vert a\Vert_{\ell^2}:=\biggl(\sum_i a_i^2\biggr)^{1/2}<\infty. $$  

The reason for excluding an exceptional curve family $\Ga_0$ in the definition
of admissibility is to guarantee that for some relevant  families of
paths $\Ga$
an admissible mass distribution exists and we have  $0<{\rm mod}_S(\Gamma)<\infty$. 

It is straightforward to check that the carpet  modulus is monotone and countably subadditive. An crucial  property of carpet modulus is its invariance under quasiconformal maps.

\begin{lemma}\label{L:Invmod}
Let $D$ be a region in $\Sph^2$, let $S$ be a carpet contained in $D$, and $\Gamma$ be a path family such that $\gamma\subseteq  D$ for each $\gamma\in\Gamma$. If  $f\co D\to\widetilde D$ is  a quasiconformal map onto another region $\widetilde D \subseteq \Sph^2$ , $\widetilde S:=f(S)$, and  $\widetilde\Gamma:=f(\Gamma)$, then 
$$
{\rm mod}_{\widetilde S}(\widetilde\Gamma)={\rm mod}_S(\Gamma).
$$
\end{lemma}

\noindent
\emph{Proof.} Note that $\widetilde S$ is also a carpet. 
Then the  peripheral circles of $S$ and of $\widetilde S$ correspond to each other under the map $f$.  So if $C_i$, $i\in \N$, is the family of peripheral circles of $S$, then $f(C_i)$, $i\in \N$, 
is the family of peripheral circles of $\widetilde S$. 

Let $\rho$ be an admissible mass distribution for $\Mod_S(\Gamma)$ with an exceptional path family $\Gamma_0$. 
Then the function $\tilde\rho$ that takes the value $\rho(C_i)$ at the peripheral circle $f(C_i)$ of $\widetilde S=f(S)$ is admissible for $\widetilde\Gamma$ with an exceptional path family $\widetilde\Gamma_0=f(\Gamma_0)$. Indeed, the $\tilde \rho$-length of every path $\tilde\gamma=f\circ \gamma$, $\gamma\in \Gamma$,  is the same as the $\rho$-length  of $\gamma$,  and the vanishing of the  conformal modulus of $\widetilde\Gamma_0$ is guaranteed by~(\ref{E:Qcqinv}).  The mass distributions  $\tilde\rho$ and $\rho$ have the same total mass.  Therefore ${\rm mod}_{\widetilde S}(\widetilde\Gamma)\leq {\rm mod}_S(\Gamma)$. We also have  the converse inequality, since $f^{-1}$ is also quasiconformal~\cite[Section~3.2, p.~17]{LV}.
\qed\medskip 



An admissible  mass distribution $\rho$ is called \emph{extremal} for ${\rm mod}_S(\Ga)$  if 
$$
{\rm mass}(\rho)={\rm mod}_S(\Ga).
$$ 
An elementary convexity argument shows that if ${\rm mod}_S(\Ga)<\infty$ and an extremal mass distribution exists, then it is unique. 
Proposition~\ref{P:infmin} below guarantees existence  of an extremal mass distribution. To prove this proposition, we need some auxiliary results. We first set up some notation.

 We let $L^2$ be the space of all  functions  $f$ on $\Sph^2$ that are square-integrable with respect to spherical measure $\sigma$, and set
$$\Vert f\Vert_{L^2}:=\biggl(\int f^2\, d\sigma\biggr)^{1/2}. $$
For two quantities $A$ and $B$ we write $A\lesssim B$ if there exists a constant $C\geq 0$ (depending on some obvious ambient parameters) such that $A\leq CB$. 

 A version of the following lemma can be found in~\cite{bB87}, see also~\cite[Exercise~2.10]{jH01}.
\begin{lemma}\label{L:Boj}  
Let  $\lambda\ge1$, and $I$ be a countable index set.  Suppose that $B_i$, $i\in I$,  is a  collection of spherical disks in $\Sph^2$, and that $a_i$, $i\in I$, are  non-negative real numbers.  Then there exists a constant $C\geq0$ that depends only on $\lambda$ such that
\begin{equation}\label{E:Boj}
\biggl\Vert \sum_{i\in I}a_i\chi_{\lambda B_i}\biggr\Vert_{L^2} \leq C\biggl \Vert \sum_{i\in I}a_i \chi_{B_i}\biggr\Vert_{L^2}.
\end{equation}
\end{lemma}

\smallskip\noindent
\emph{Proof.} We may assume that $\sum_{i\in I}a_i \chi_{B_i}\in L^2$. 
Let $\phi\in L^2$,  and $M(\phi)$ denote the uncentered maximal function of $\phi$ (for the  definition and the basic properties of the maximal function operator see
 \cite[Chapter 1]{eS70}). 
Then there is an absolute constant $c$ such that
$$
\begin{aligned}
&\biggl|\int\biggl(\sum_{i\in I} a_i\chi_{\lambda B_i}\biggr)\phi\, d\sigma\biggr|
=\biggl|\sum_{i\in I}a_i\int_{\lambda B_i}\phi\, d\sigma\biggr|\\
&\leq c\lambda^2\sum_{i\in I}a_i\int_{B_i}M(\phi)\,d\sigma
=c\lambda^2\int\biggl(\sum_{i\in I} a_i\chi_{B_i}\biggr)M(\phi)\, d\sigma\\
&\leq c \lambda^2\biggl\Vert \sum_{i\in I}a_i\chi_{B_i}
\biggr\Vert_{L^2}  \Vert M(\phi)\Vert_{L^2}.
\end{aligned}
$$
It is known  (see e.g.,~\cite[Theorem 1(c), p.~5]{eS70})
that the maximal function satisfies the inequality
$$
\Vert M(\phi)\Vert_{L^2}\leq H\Vert\phi\Vert_{L^2},
$$
where $H$ is an absolute constant. The self-duality of $L^2$ now gives
inequality~(\ref{E:Boj}) with $C=cH\lambda^2$.
\qed\medskip 

A \emph{quasicircle} in a metric space $X$ is a Jordan curve that is quasisymmetrically equivalent  to the unit circle in $\R^2$ equipped with the Euclidean metric. 
We say that  a family $\{C_i:i\in I\}$ of Jordan curves in a metric space $X$ consists of  \emph{uniform quasicircles},  
if
there exists a homeomorphism $\eta\co [0,\infty)\to[0,\infty)$ such that every curve  $C_i$ in the  family is 
the image of  an $\eta$-quasisymmetric map of the unit circle.

\begin{lemma}\label{lem:modzero}  
Let $S$ be a carpet in $\Sph^2$ whose peripheral circles are uniform quasicircles, and   let  $\Gamma$ be   a path family in $\Sph^2$.
If    ${\rm mod}_S(\Gamma)=0$,
then  ${\rm mod}(\Gamma)=0$. 
\end{lemma}

\smallskip\noindent
\emph{Proof.} We write   $S$ as in \eqref{eq:carprep}, and set $C_i=\partial D_i$ for $i\in \N$. Under the given hypotheses suppose that  ${\rm mod}_S(\Gamma)=0$. It follows from the definitions that then $\Gamma$ cannot contain any constant paths. 
 So if  for $k\in \N$ we define 
$$
\Gamma_{k}:=\{\gamma\in\Gamma\co {\rm diam}(\gamma)\geq1/k\}, 
$$
then $\Gamma=\bigcup_{k}\Gamma_k$, and it is enough to show that ${\rm mod}(\Gamma_k)=0$ for every $k\in\N$.   By monotonicity of carpet modulus we have ${\rm mod}_S(\Gamma_k)=0$. This means  that we are actually reduced to proving the statement of the lemma under the additional 
assumption that ${\rm diam}(\gamma)\geq\delta$ for all $\gamma\in \Gamma$, where $\delta>0$. 

Since ${\rm mod}_S(\Gamma)=0$, for each $n\in \N$ there exists 
$\rho_n\in \ell^2$ with $\Vert \rho_n\Vert_{\ell^2}<1/2^n$
 and an exceptional family $\widetilde \Gamma_n\sub 
\Gamma$ with  ${\rm mod}(\widetilde \Gamma_n)=0$ such that  
$$\sum_{\ga\cap C_i\neq\emptyset}\rho_n(C_i)\geq 1$$
for all $\gamma\in \Gamma\setminus \widetilde \Gamma_n. $
Let $\rho=\sum_n\rho_n$. Then 
$$ {\rm mass}(\rho)=\Vert \rho\Vert_{\ell^2}^2<\infty. $$
Moreover, if $\widetilde \Gamma:= \bigcup_n \widetilde \Gamma_n$,
then ${\rm mod}( \widetilde \Gamma)=0$ and  
\begin{equation}\label{eq:totalinf}
\sum_{\ga\cap C_i\neq\emptyset}\rho(C_i)=\infty
\end{equation} 
for all $\gamma\in \Gamma\setminus \widetilde \Gamma.$

 Since we assume that the peripheral circles $C_i=\partial D_i$ of $S$ are uniform quasicircles, there exists $\lambda\geq 1$ with the following property (see, e.g., \cite[Proposition~4.3]{mB04}):  
for each $i\in \N$ there 
exists $x_i\in\Sph^2$,  and $0<r_i\leq R_i$ with $R_i/r_i\le \lambda$ such that 
$$
 B(x_i,r_i)\subseteq D_i\subseteq B(x_i,R_i). 
$$

Now we consider the  density $\tilde \rho$ on the sphere defined by
$$
\tilde \rho=\sum_{i}\frac{\rho(C_i)}{R_i}\chi_{B(x_i,2R_i)}. 
$$
If $4R_i<\delta$,  then  every path $\gamma\in\Gamma$ that meets  $C_i$ 
must also meet the complement of $B(x_i,2R_i)$, since ${\rm diam}(\gamma)\ge \delta$. So $\gamma$ 
will meet  both  complementary components of $B(x_i,2R_i)\setminus\overline B(x_i,R_i)$. If $\gamma$ is locally rectifiable, this 
implies 
$$ \int_\gamma \chi_{B(x_i,2R_i)}\, ds\ge R_i. $$
Since there are only finitely many $i\in \N$ with $4R_i>\delta$, it follows  
 from~\eqref{eq:totalinf} that 
$$
\int_{\gamma}\tilde \rho\,ds=\infty
$$
for every locally rectifiable path $\gamma\in \Gamma\setminus \widetilde \Gamma$.
On the other hand, by Lemma~\ref{L:Boj},
$$
\int\tilde \rho^2\, d\sigma\lesssim \sum_{i}\frac{\rho(C_i)^2}{R_i^2}\sigma(B(x_i,r_i))\lesssim {\rm mass}(\rho)<\infty.
$$
 This implies ${\rm mod}(\Gamma\setminus \widetilde \Gamma) = 0$, and so
 $${\rm mod}(\Gamma)\le {\rm mod}(\Gamma\setminus \widetilde \Gamma)
 +{\rm mod}(\widetilde \Gamma)=0. $$
 Hence ${\rm mod}(\Gamma)=0$ as desired.  \qed \medskip

\begin{proposition}\label{P:infmin}
Let $S$ be a carpet in $\Sph^2$ whose peripheral circles are uniform quasicircles,  and  let $\Gamma$ be an arbitrary path family in $\Sph^2$ with ${\rm mod}_S(\Gamma)<\infty$. Then the extremal mass distribution for $
{\rm mod}_S(\Gamma)$ exists, i.e., the infimum in the definition of ${\rm mod}_S(\Gamma)$ is attained as a minimum. 
\end{proposition}
\smallskip\noindent
\emph{Proof.} 
Let $(\rho_n)$ be a sequence of admissible mass distributions for ${\rm mod}_S(\Gamma)$ such that ${\rm mass}(\rho_n)\to{\rm mod}_S(\Gamma)$ as $n\to\infty$. 
If $C_i$, 
$i\in \N$, are the peripheral circles of $S$, then each $\rho_n$ is given by the  sequence $\rho_n=(\rho_n(C_i))$  of weights it assigns to the peripheral 
circles. 

Since  ${\rm mod}_S(\Gamma)<\infty$, we may assume that  for a suitable constant we  have  ${\rm mass}(\rho_n)\leq C$ for all $n$. By passing to a subsequence using a  standard diagonalization argument, we may also assume that the limit 
$$
\rho(C_i):=\lim_{n\to\infty}\rho_n(C_i)
$$ exists for each  $i\in \N$. 
We claim that the mass distribution $\rho=(\rho(C_i))$ is extremal.  

First, it is clear that ${\rm mass}(\rho)\leq{\rm mod}_S(\Gamma)$. Indeed, for every $\epsilon>0$ and $m\in \N$ there exists $N\in\N$ such that 
$$
\sum_{i=1}^m\rho_n(C_i)^2\leq{\rm mod}_S(\Gamma)+\epsilon
$$  for all $n\geq N$.
Taking the limit as $n\to\infty$, we get 
$$
\sum_{i=1}^m\rho(C_i)^2\leq {\rm mod}_S(\Gamma)+\epsilon
$$  for all $m\in\N$.
Since $m$ and $\epsilon$ are arbitrary, this gives ${\rm mass}(\rho)\leq {\rm mod}_S(\Gamma)$.

To complete the proof we have  to show that $\rho$ is admissible, which is harder to establish. By Mazur's Lemma (see, e.g.,~\cite[Theorem~2, p.~120]{kY80}) 
 there is a sequence of convex combinations $(\tilde\rho_N)$, where 
$$
\tilde\rho_N=\sum_{n=1}^N\lambda_n^N\rho_n,\quad \lambda_n^N\geq0,\quad \sum_{n=1}^N\lambda_n^N=1,
$$ 
that converges to $\rho$ in $\ell^2$. Every element of the sequence $(\tilde\rho_N)$ is admissible for $\Gamma$, where the exceptional path family $\widetilde\Ga_N$ for $\tilde\rho_N$ is the union of the exceptional path families for $\rho_n,\ n=1,2,\dots, N$. Since $(\tilde\rho_N)$ converges to $\rho$ in $\ell^2$,  it is also a minimizing sequence for ${\rm mod}_S(\Gamma)$. 

By passing to a subsequence, we may assume that 
\begin{equation}\label{E:Geom}
||\tilde\rho_{N}-\rho||_{\ell^2}\leq \frac1{2^N},\quad \end{equation}
for all $N\in\N$.

Let 
$$
\Gamma_\infty=\biggl\{\gamma\in\Gamma\co \limsup_{N\to\infty}\sum_{\gamma\cap C_i\neq\emptyset}|\tilde\rho_N(C_i)-\rho(C_i)|\neq0\biggr\}
$$
and
$$
\Gamma_N=\bigg\{\gamma\in\Gamma\co \sum_{\gamma\cap C_i\neq\emptyset}|\tilde\rho_N(C_i)-\rho(C_i)|\geq\frac1N\bigg\}.
$$
Then $\Gamma_\infty\subseteq\bigcap_{n}\bigcup_{N\ge n}\Gamma_N$. 
Indeed, let $\gamma\in\Ga_\infty$ be arbitrary. Then there exists $\delta>0$ and a sequence of natural numbers $(N_k)$ with $N_k\to \infty$ as $k\to \infty$  such that 
$$
\sum_{\ga\cap C_i\neq\emptyset}|\tilde\rho_{N_k}(C_i)-\rho(C_i)|\geq\delta
$$
for all $k\in\N$. Now let $n$ be arbitrary. We choose $k$ so large that $N_k\geq n$ and $1/N_k\leq\delta$. Then $\ga\in \Ga_{N_k}
\sub \bigcup_{N\ge n}\Gamma_N$. Hence $\ga\in \bigcap_n \bigcup_{N\ge n}\Gamma_N$ as desired. 

It follows that the mass distributions 
$$
\rho_{\infty,n}=\sum_{N=n}^\infty N|\tilde\rho_N-\rho|,\ n=1,2,\dots, 
$$
are admissible for $\Mod_S(\Gamma_\infty)$. Since ${\rm mass}(\rho_{\infty,n})\to0$ as $n\to\infty$ by~(\ref{E:Geom}), this implies that   ${\rm mod}_S(\Gamma_\infty)=0$.
Invoking  Lemma~\ref{lem:modzero} we conclude  that ${\rm mod}(\Gamma_\infty)=0$. 

If $\gamma$ is in $\Gamma\setminus(\Gamma_\infty\cup \bigcup_{N}\widetilde\Ga_N)$, then 
$$
\begin{aligned}
\sum_{\gamma\cap C_i\neq\emptyset}\rho(C_i)&\geq\limsup_{N\to\infty} \biggl( \sum_{\gamma\cap C_i\neq\emptyset}\tilde\rho_N(C_i)-\sum_{\gamma\cap C_i\neq\emptyset}|\tilde\rho_N(C_i)-\rho(C_i)|\biggr)\\
&\geq1-\limsup_{N\to\infty} \biggl(\sum_{\gamma\cap C_i\neq\emptyset}|\tilde\rho_N(C_i)-\rho(C_i)|\biggr)=1.\end{aligned}
$$
Moreover,
$${\rm mod}\biggl (\Gamma_\infty\cup \bigcup_{N}\widetilde\Ga_N\biggr)\le {\rm mod}(\Gamma_\infty)
+\sum_N  {\rm mod}(\widetilde \Gamma_N)=0.$$ 
It follows that $\rho$ is admissible for $\Gamma$ as desired.  
\qed\medskip 

\section{Carpet modulus with respect to a group}\label{S:Cmodgp}

\no 
Let $S$ be a carpet in $\Sph^2$. In this section we assume that $S$ is written as in \eqref{eq:carprep}, and that $C_i=\partial D_i$,  $i\in \N$, denotes the peripheral circles of $S$. Let 
$G$ be a group  of homeomorphisms   of $S$. If $g\in G$ and $C\sub S$ is a peripheral circle of $S$, then $g(C)$ is also a peripheral circle of 
$S$. So the whole orbit  $\mathcal {O}=\{g(C):g\in G\}$ 
of $C$ under the action of $G$ consists of peripheral circles 
of $S$.  If $\Gamma$ is a family of paths in $\Sph^2$, we define the  carpet modulus ${\rm mod}_{S/G}(\Gamma)$ of $\Gamma$ with respect to the action of $G$ as follows. 
A \emph{(invariant) mass distribution} $\rho$ is a non-negative function defined on the peripheral circles of $S$ that takes the same value  on each  peripheral circle in the same orbit;  so $\rho(g(C))=\rho(C)$ for all $g\in G$ and all peripheral circles $C$ of $S$.  
 Such a  mass distribution is \emph{admissible} for ${\rm mod}_{S/G}(\Gamma)$
if there exists an {\em exceptional family} $\Gamma_0\sub \Gamma$ with ${\rm mod}(\Gamma_0)=0$ and 
\begin{equation}\label{eq:admodgroup}
\sum_{\gamma\cap C_i\neq\emptyset}\rho(C_i)\geq1\end{equation} 
for all  $\gamma\in\Gamma\setminus\Gamma_0$.

If $\rho$ is a mass distribution and $\mathcal{O}$  an  orbit of  peripheral circles, we set 
$\rho(\mathcal{O}):=\rho(C)$, where $C\in \mathcal{O}$. We  define the \emph{(total) mass} of $\rho$  as
$$
{\rm mass}_{S/G}(\rho)=\sum_{\mathcal{O}}\rho(\mathcal{O})^2,
$$
where the sum is taken over all orbits of peripheral circles 
under the action of $G$.


The \emph{carpet modulus of} $\Gamma$ \emph{with respect to the group} $G$
is defined as 
$$
{\rm mod}_{S/G}(\Gamma)=\inf_\rho \,  \{{\rm mass}_{S/G}(\rho)\},$$
where the infimum is taken over all admissible mass distributions
$\rho$.  An admissible mass distribution $\rho$ realizing this infimum is called {\em extremal} for ${\rm mod}_{S/G}(\Gamma)$. 

Note that each orbit contributes with exactly one term  to the total mass of a mass distribution. In contrast,   the admissibility condition is similar to the one  for  carpet modulus: each 
 peripheral circle that intersects the curve $\ga$ contributes a  term  to the sum in \eqref{eq:admodgroup} and we may get multiple contributions from each orbit;  we restrict ourselves though to invariant mass distributions that are constant 
 on each orbit of the action of $G$ on peripheral circles. 
 
 Carpet modulus with respect to a group 
 has similar monotonicity and subadditivity properties as  carpet modulus and conformal modulus.  We formulate its invariance property under quasiconformal maps explicitly.

\begin{lemma}\label{L:Invgpmod}
Let $D$ be a region in $\Sph^2$,  $S$ a carpet contained in $D$, $G$  a group of homeomorphisms  on $S$,  and $\Gamma$  a path family such that $\gamma\subseteq  D$ for each $\gamma\in\Gamma$. 
If  $f\co D\to\widetilde D$ is  a quasiconformal map onto another region 
$\widetilde D\subseteq \Sph^2$, 
 and we define   $\widetilde S:=f(S)$,   $\widetilde\Gamma:=f(\Gamma)$, and 
$\widetilde G:=(f|S)\circ G\circ (f|S)^{-1}, $ then 
$$
{\rm mod}_{\widetilde S/\widetilde G}(\widetilde \Gamma)={\rm mod}_{S/G}(\Gamma).
$$
\end{lemma}

\noindent
\emph{Proof.} Note that $\widetilde S$ is a carpet, and $f|S$ is 
homeomorphism from $S$ onto $\widetilde S$. Hence 
$\widetilde G$ is a group of homeomorphisms on $\widetilde S$. 
The argument is now along the same lines as the proof of   Lemma~\ref{L:Invmod}
and we omit the details. 
\qed\medskip 

The following proposition gives a criterion for the existence of an extremal mass distribution for carpet modulus with respect to the group. 

\begin{proposition}\label{P:infmingroup}
Let $S$ be a carpet in $\Sph^2$ whose peripheral circles are uniform quasicircles,  let  $G$ be a group of homeomorphisms of $S$, and let $\Gamma$ be a path family in $\Sph^2$ with ${\rm mod}_{S/G}(\Gamma)<\infty$.

Suppose that  for each $k\in \N$ there exists a family of peripheral circles $\mathcal{C}_k$ of $S$ and a constant $N_k\in\N$
with the following properties:

\begin{itemize} 

\smallskip
\item[{(i)}] if $\mathcal{O}$ is any orbit of peripheral circles of 
$S$ under the action of $G$, then $\#(\mathcal{O}\cap 
\mathcal{C}_k)\le N_k$ for all $k\in \N$, 

\smallskip 
\item[{(ii)}] if $\Gamma_k$ is the family of all paths in $\Gamma$ that only meet peripheral circles in $\mathcal{C}_k$, then 
$\Gamma=\bigcup_k \Gamma_k$. \end{itemize} 

 Then an  extremal mass distribution for ${\rm mod}_{S/G}(\Gamma)$ exists, i.e., the infimum in the definition of ${\rm mod}_{S/G}(\Gamma)$ is attained. 
\end{proposition}

\no
\emph{Proof.}  We first observe that  one has an analog of Lemma~\ref{lem:modzero}; namely, if in addition to  the given hypotheses  ${\rm mod}_{S/G}(\Gamma)=0$, 
then  ${\rm mod}(\Gamma)=0$.  

The proof of this implication is very similar to the proof of 
Lemma~\ref{lem:modzero}. 
As in the proof of this lemma, we can  make the additional assumption that there exists $\delta>0$ such that $\diam(\gamma)\ge \delta$ for all $\gamma\in \Gamma$. 
Using  that 
${\rm mod}_{S/G}(\Gamma)=0$ one can find an invariant  mass distribution  $\rho$ 
with ${\rm  mass}_{S/G}(\rho)<\infty$ and a family $\widetilde \Gamma\sub \Gamma$ with $\Mod(\widetilde\Gamma)=0$ such that 
\begin{equation}\label{eq:suminfty1}
\sum_{\ga\cap C_i\neq\emptyset}\rho(C_i)=\infty\end{equation} 
for all $\ga\in\Ga\setminus\widetilde\Ga$. 

 Since peripheral circles of $S$, represented as in \eqref{eq:carprep},  are uniform quasicircles, there exists $\lambda\geq 1$ such that
for each $i\in \N$ we can find 
 $x_i\in\Sph^2$ and $0<r_i\leq R_i$ with 
$$
 B(x_i,r_i)\subseteq D_i\subseteq B(x_i,R_i),
$$
and 
$R_i/r_i\leq\lambda$ . 

Now fix $k\in \N$,  consider the family $\Gamma'_k:=(\Gamma\setminus\widetilde\Ga)\cap \Gamma_k$ of all paths in $\Ga\setminus\widetilde\Ga$ that only intersect peripheral circles in the family
$\mathcal{C}_k$, and  
let
$$
\tilde \rho=\sum_{C_i\in\mathcal{C}_k}\frac{\rho(C_i)}{R_i}\chi_{B(x_i,2R_i)}. 
$$
Using our   hypothesis (i) and Lemma~\ref{L:Boj} we see that 
$$\int \tilde \rho^{2}\,d\sigma\lesssim \sum_{C_i\in\mathcal{C}_k}\rho(C_i)^2\le N_k {\rm mass}_{S/G}(\rho)<\infty.$$
On the other hand, similarly as in the proof of  Lemma~\ref{lem:modzero},   by  \eqref{eq:suminfty1} we have 
$$\int_\gamma\tilde \rho\, ds =\infty$$ 
for every locally rectifiable path $\ga \in \Gamma'_k$. It follows that ${\rm mod}(\Gamma_k')=0$. 

Our hypothesis (ii) implies that 
$\Gamma=\widetilde \Gamma\cup \bigcup_k\Gamma'_k$. Since 
all families in the last union have vanishing modulus, we conclude 
${\rm mod}(\Gamma)=0$ as desired.

 Now  the proof of the statement is almost identical  to the proof of  Proposition~\ref{P:infmin}. The only difference is that we use mass distributions that are constant on each orbit of peripheral circles under the action of $G$, and  that for control on the masses   of the relevant distributions we use an    $\ell^2$-space  indexed by these orbits.
  Note that hypothesis (ii) passes to every subfamily 
 of $\Gamma$, so we can apply  the first part of the proof to the family that corresponds to $\Gamma_\infty$ in the proof of 
 Proposition~\ref{P:infmin}. We omit the details. 
\qed\medskip

If $\psi$ is a homeomorphism of  the carpet $S\sub \Sph^2$, we denote by $\langle\psi\rangle$ the cyclic group of homeomorphisms on $S$ generated by $\psi$.  If $\Gamma$ is a path family in $\Sph^2$ and $\Psi$ is a homeomorphism on $\Sph^2$, then  $\Gamma$ is called $\Psi$-{\em invariant} if $\Psi(\Gamma)=\Gamma$.  
The following lemma  gives a precise relationship between the carpet modulus with respect to a cyclic group and its subgroups. 

\begin{lemma}\label{L:Cyclgp}
Let $S$ be a carpet in $\Sph^2$, let $\Psi\: \Sph^2\ra \Sph^2$ be a quasiconformal map with $\Psi(S)=S$, and define $\psi:=\Psi|S$.
Suppose that  $\Gamma$ is  a $\Psi$-invariant path  family in $\Sph^2$ such that  for every peripheral circle $C$ of $S$ that meets some  path in $\Gamma$  we have
$
\psi^n(C)\neq C
$ 
for all $n\in \Z\setminus \{0\}$. 

Then
$$
{\rm mod}_{S/\langle\psi^k\rangle}(\Gamma)=k\,{\rm mod}_{S/\langle\psi\rangle}(\Gamma)
$$ for every $k\in \N$. 
\end{lemma}
\smallskip\noindent

\no
\emph{Proof.} Fix $k\in \N$.  Let $\epsilon>0$,  and $\rho$ be an admissible mass distribution for ${\rm mod}_{S/\langle\psi\rangle}(\Gamma)$  such that 
$$
{\rm mass}_{S/\langle\psi\rangle}(\rho)\leq {\rm mod}_{S/\langle\psi\rangle}(\Gamma)+\epsilon.
$$
Here it is understood that $\rho$ is invariant in the sense that it is constant on orbits of peripheral circles under the action of $\langle\psi\rangle$.
Then $\rho$ is also constant on  orbits under the action of 
$\langle\psi^k\rangle$, and hence admissible for  
${\rm mod}_{S/\langle\psi^k\rangle}(\Gamma)$. 

Each orbit of a peripheral circle  under the action of 
$\langle\psi\rangle$ consists of at most $k$ orbits under the action of $\langle\psi^k\rangle$. Therefore, we obtain
\begin{eqnarray*} 
{\rm mod}_{S/\langle\psi^k\rangle}(\Gamma)&\leq & {\rm mass}_{S/\langle\psi^k\rangle}(\rho)\\ &\le& k\,{\rm mass}_{S/\langle\psi\rangle}(\rho)\leq k\,{\rm mod}_{S/\langle\psi\rangle}(\Gamma)+k\epsilon.
\end{eqnarray*} 
Since $\epsilon$ is arbitrary, we conclude
\begin{equation*}
{\rm mod}_{S/\langle\psi^k\rangle}(\Gamma)\leq k\,{\rm mod}_{S/\langle\psi\rangle}(\Gamma).
\end{equation*}

Conversely, let $\epsilon>0$,  and $\rho$ be an admissible  invariant mass distribution for ${\rm mod}_{S/\langle\psi^k\rangle}(\Gamma)$ such that
$$
{\rm mass}_{S/\langle\psi^k\rangle}(\rho)\leq {\rm mod}_{S/\langle\psi^k\rangle}(\Gamma)+\epsilon.
$$

Note that if $C$ is a peripheral circle of $S$ and 
no path in $\Gamma$ meets $C$, then by $\Psi$-invariance of 
$\Gamma$ no path in $\Gamma$ meets any of  the peripheral circles in the orbit of $C$ under $\langle\psi\rangle$. 
This implies that we may assume that $\rho(C)=0$ for all peripheral circles $C$ that do not meet any path in $\Gamma$. 
 
 Consider $\tilde\rho$ given by 
$$
\tilde\rho=\frac 1k(\rho+\rho\circ\psi+\dots+\rho\circ\psi^{k-1}).
$$
Here $\rho\circ \psi^j$ denotes the  mass distribution that assigns  the value $\rho(\psi^j(C))$ to  a  peripheral $C$ of $S$.
 
We have that $\rho\circ \psi^k=\rho$, since $\rho$ is constant on orbits of peripheral circles under the action of $\psi^k$. 
This implies that $\tilde \rho\circ \psi=\tilde \rho$ and so $\tilde\rho$ 
is constant on orbits of peripheral circles under the action of $\psi$. 

Let  $\Gamma_0$ be  an exceptional family for $\rho$, 
 and define 
$$\widetilde\Gamma_0=\bigcup_{n\in \{ -(k-1), \dots, -1,0\} }\Psi^n(\Gamma_0). 
$$
Since $\Mod(\Gamma_0)=0$, we  have 
$\Mod(\widetilde\Gamma_0)=0$. Moreover, the $\Psi$-invariance of $\Gamma$ implies  that 
$$ \sum_{\gamma\cap C_i\ne \emptyset} \tilde \rho(C_i) \ge 1 $$ 
for all $\gamma\in \Gamma\setminus \widetilde\Gamma_0$.
Hence $\tilde \rho$ is  admissible for ${\rm mod}_{S/\langle\psi\rangle}(\Gamma)$. 

It follows that 
$$
{\rm mod}_{S/\langle\psi\rangle}(\Gamma)\leq{\rm mass}_{S/\langle\psi\rangle}(\tilde\rho). 
$$
Since $\rho$ assigns $0$ to all peripheral circles of $S$ that 
do not meet any path in $\Gamma$, the same is true for $\tilde \rho$. So if $C$ is a peripheral circle of $S$ with $\tilde \rho(C)\ne 0$, then $C$ meets some  path in $\Gamma$ and so  our hypotheses imply that the peripheral circles $\psi^n(C)$, $n\in \Z$, are all distinct. Hence the $\langle\psi\rangle$-orbit of $C$ consists of precisely $k$ orbits of $C$ under $\langle\psi^k\rangle$. It follows 
that 
$$ {\rm mass}_{S/\langle\psi^k\rangle}(\tilde\rho)=k\, {\rm mass}_{S/\langle\psi\rangle}(\tilde\rho). $$ 
Moreover,  the convexity of the norm in $\ell^2$ implies that 
 $$
{\rm mass}_{S/\langle\psi^k\rangle}(\tilde\rho)\leq{\rm mass}_{S/\langle\psi^k\rangle}(\rho).
$$
We conclude 
\begin{eqnarray*}
k\, {\rm mod}_{S/\langle\psi\rangle}(\Gamma)&\leq&k\, {\rm mass}_{S/\langle\psi\rangle}(\tilde\rho)\,=\,{\rm mass}_{S/\langle\psi^k\rangle}(\tilde\rho)\\ &\leq& {\rm mass}_{S/\langle\psi^k\rangle}(\rho)\, \le\, {\rm mod}_{S/\langle\psi^k\rangle}(\Gamma)+{\epsilon}.
\end{eqnarray*}
Since $\epsilon$ was  arbitrary, this gives the other desired inequality 
$$k\,{\rm mod}_{S/\langle\psi\rangle}(\Gamma)\le {\rm mod}_{S/\langle\psi^k\rangle}(\Gamma).$$ The statement follows. \qed\medskip 

%
%

\section{Auxiliary results}\label{S:Aux} 

\no
In this section  we state  results from~\cite{mB04} and \cite{BKM06} that are used in this paper and derive some consequences. 

\begin{proposition}\label{P:Qcext}
Let $S$ be a carpet in $\Sph^2$ whose peripheral circles are uniform quasicircles and let $f$ be a quasisymmetric map of $S$ onto another  carpet $\widetilde S\subseteq\Sph^2$. Then there exists a quasiconformal map $F$ on $\Sph^2$ whose restriction to $S$ is $f$.
\end{proposition}

This follows immediately from~\cite[Proposition~5.1]{mB04}.


Suppose  $\{C_i:i\in I\}$ is a family of continua  in a metric space $X$, i.e., each set $C_i$ is a compact connected set consisting of more than one point.  These continua are said to be \emph{uniformly relatively 
separated} if the pairwise relative distances are uniformly bounded away from zero, i.e., there exists $\delta>0$ such that 
$$
\Delta(C_i, C_j):=\frac{{\rm{dist}}(C_i,
C_j)}{\min\{{\rm{diam}}(C_i),{\rm{diam}}(C_j)\}}\geq\delta
$$
for any two distinct $i$ and $j$. 
The uniform relative separation property is preserved  under quasisymmetric maps, see~\cite[Corollary~4.6]{mB04}. 


Recall that a carpet $S\sub \Sph^2$ is called {\em round} if its peripheral circles are geometric circles. So if $S$ is written as in \eqref{eq:carprep}, then
each Jordan region $D_i$ is an open spherical disk. 

\medskip 
\noindent  \begin{theorem}[Uniformization by round carpets]\label{T:rdcarp}
If $S$ is a carpet in $\Sph^2$ whose peripheral circles are uniformly relatively  separated uniform quasicircles, 
then there exists a quasisymmetric map of $S$ onto a round carpet. 
\end{theorem}
\medskip 

This is  \cite[Corollary~1.2]{mB04}.

\medskip 
\noindent \begin{theorem}[Quasisymmetric rigidity of round carpets]\label{T:Riground}  Let  $S$ be  a round carpet in $\Sph^2$ of measure zero.  
Then every quasisymmetric map of 
$S$ onto any other round carpet is the restriction of a M\"obius transformation. 
\end{theorem}

\medskip 
This is
~\cite[Theorem~1.2]{BKM06}. Here by definition a {\em M\"obius 
transformation} is a fractional linear transformation on $\Sph^2\cong \widehat \C$, or the complex-conjugate of such a map. So we allow a 
M\"obius 
transformation to be orientation-reversing. 

Let $S\sub \Sph^2$ be a carpet, and $f\: S\ra \Sph^2$ be a homeomorphic embedding. Then $f$ has a homeomorphic extension to a homeomorphism $F\: \Sph^2\ra \Sph^2$ (see the proof of Lemma~5.3 in \cite{mB04}). We call $f$ {\em orientation-preserving} if $F$ (and hence every homeomorphic extension of $f$) is orientation-preserving on $\Sph^2$.  On a more intuitive level, the map $f$ is orientation-preserving if the following is true: if we orient any    peripheral circle $C$ of $S$ so that $S$ lies ``to the left" of $C$ with this orientation, then the image carpet 
$f(S)$ lies ``to the left" of its  peripheral circle $f(C)$ equipped with the induced orientation.

\begin{corollary}[Three-Circle Theorem]\label{C:TCT}
Let $S$ be a carpet in $\Sph^2$ of  measure zero whose peripheral circles are uniformly relatively separated uniform quasicircles. Let $C_1, C_2,C_3$ be three distinct peripheral circles of $S$. 
If $f$ and $g$ are two orientation-preserving 
quasisymmetric self-maps of $S$ such that $f(C_i)=g(C_i)$ for $i=1,2,3$, 
then $f=g$.
\end{corollary}

This follows from  \cite[Theorem~1.5]{mB04} applied to $f^{-1}\circ g$. 

\begin{corollary}\label{C:Cycl} Let $S$ be a carpet in $\Sph^2$ of  measure zero whose peripheral circles are uniformly relatively separated uniform quasicircles. Let $C$ be a peripheral circle of $S$, and $p,q$ two distinct points on $C$, and  $G$ be the  group  of all orientation-preserving quasisymmetric self-maps of $S$ that fix $p$ and $q$. Then 
$G=\{\id_S\}$ or $G$ is an infinite cyclic group. 
 \end{corollary}

In other words, either $G$ is trivial or isomorphic to $\Z$. 

\medskip 
\noindent
\emph{Proof.}
By Theorem~\ref{T:rdcarp} there exists  a quasisymmetric map $f$ of $S$ onto a round carpet $\widetilde S$. 
Using Proposition~\ref{P:Qcext} we can extend  $f$  to a quasiconformal map on $\Sph^2$. Since quasiconformal maps preserve the class of sets of measure 
zero~\cite[Theorem~1.3, p.~165]{LV},
 $\widetilde S$ has measure zero as well.
According to Theorem~\ref{T:Riground}, 
the conjugate group $\widetilde G=f\circ G\circ f^{-1}$
consists of the restrictions of orientation-preserving M\"obius transformations with two fixed points $\tilde p,\tilde q$ on a peripheral circle $\widetilde C$ of $\widetilde S$.  
By post-composing $f$ with a M\"obius transformation we may assume that  $\tilde p=0$,  $\tilde q=\infty$, and that $\widetilde C$ is the extended real line.  Moreover, we may assume that $\widetilde S$ is contained in the upper half-plane. Then the  maps  in $\widetilde G$ are of the form $z\mapsto \lambda z$ with $\lambda> 0$.   The multiplicative group of the factors 
$\lambda$ that arise in this way must be  discrete (this follows from the fact that peripheral circles are mapped to peripheral circles)  and hence forms 
a cyclic group. It follows that  $\widetilde G$, and hence also $G$, is the 
trivial group consisting only of the identity map or an infinite cyclic group.  
\qed\medskip

\begin{corollary}\label{C:CC} Let $S$ be a carpet in $\Sph^2$ of  measure zero whose peripheral circles are uniformly relatively separated uniform quasicircles. Let $C_1$ and $C_2$ be two distinct peripheral circle of $S$, and $G$ be the  group  of all orientation-preserving quasisymmetric self-maps of $S$ that fix $C_1$ and $C_2$ setwise. Then 
$G$ is a finite cyclic group.  \end{corollary}

\noindent
\emph{Proof.} As in the proof of Corollary~\ref{C:Cycl} we can reduce to the case that $S$ is a round carpet of measure zero. By applying an auxiliary M\"obius-transformation if necessary, we may also assume that $C_1$ and $C_2$ are Euclidean circles both centered at $0$. Then by Theorem~\ref{T:Riground} each ele\-ment in $G$ is (the restriction of) a rotation around $0$. Moreover, $G$ must be a discrete group as it maps peripheral circles of $S$ to peripheral circles. Hence $G$ is finite cyclic.   \qed\medskip 

\begin{corollary}\label{C:Ccp} Let $S$ be a carpet in $\Sph^2$ of  measure zero whose peripheral circles are uniformly relatively separated uniform quasicircles,  
 $C_1$ and $C_2$ be two distinct peripheral circles of $S$, and $p\in S$. If $f$ is an orientation-preserving quasisymmetric self-map  of $S$ such that $f(C_1)=C_1$, 
 $f(C_2)=C_2$,  and $f(p)=p$, then $f$ is the identity on $S$. 
\end{corollary}

\noindent
\emph{Proof.} By the argument as in the proof of Corollary~\ref{C:CC} we can reduce to the case where $S$ is a round 
carpet,  $C_1$ and $C_2$ are circles both centered at $0$, and $f$ is a rotation around $0$.   Since  $C_1$ and $C_2$ are distinct perpheral circles of $S$, these sets are disjoint and  bound two  disjoint Jordan regions. This implies that $S$ is contained in the Euclidean annulus with boundary components $C_1$ and $C_2$. Since $p\in S$, it follows that $p\ne 0,\infty$. Since $f$ is a rotation around $0$ and fixes $p$, the map   $f$ 
must be  the identity on $S$.  \qed\medskip

The metric on $\C^*=\C\setminus \{0\}$ induced by the length element 
$|dz|/|z|$  is called the {\em flat metric (on $\C^*$)}. Equipped with this metric, 
$\C^*$ is isometric to an infinite cyclinder of circumference $2\pi$. 
The following terminology is suggested by this geometric picture.  

A  $\C^*$-\emph{cylinder} $A$ is a set of the form
$$
A=\{z\in\C\co r\leq|z|\leq R\}, 
$$
where $0<r<R<\infty$. The boundary components $\{z\in \C: |z|=r\}$ and $\{z\in \C: |z|=R\}$ are called the \emph{inner} and \emph{outer} boundary com\-ponents of $A$, respectively.  
A $\C^*$-\emph{square}  $Q$ is a    Jordan region of the form 
$$
Q=\{\rho e^{\iu \theta}\co a<\rho<b,\ \alpha<\theta<\beta\}, 
$$
 where 
$0<\log(b/a)=\beta-\alpha<2\pi$.
We call  the quantity 
$$\ell_{\C^*}(Q):=\log(b/a)=\beta-\alpha$$ its \emph{side length}, 
the  set  $\{ a e^{\iu\theta}\co  \alpha\le \theta\le \beta \}$ the {\em bottom side},  and the set $\{ a e^{\iu\theta}\co  \alpha\le \theta\le \beta \}$ the {\em top side}  of $Q$. The sets
 $\{ \rho  e^{\iu\alpha}\co  a\le \rho\le b \}$ and $\{ \rho  e^{\iu\beta}\co  a\le \rho\le b \}$ are referred to as the {\em vertical sides} of 
 $Q$. The {\em corners} of $Q$ are the four points that are endpoints of one of the sides  of $Q$. 
 
A \emph{square carpet} $T$ {\em in a $\C^*$-cylinder $A$} is a carpet that can be written as 
$$
T=A\setminus\bigcup_iQ_i,
$$ 
where the sets  $Q_i$, $i\in I$, are  $\C^*$-{squares} whose closures are pairwise disjoint and  contained in the interior of  $A$.  Very similar  terminology was  employed in \cite{mB04}. Note that in contrast to  \cite{mB04}  our 
$\C^*$-cylinders  $A$ are closed   and the  $\C^*$-squares $Q$ are open sets.  

\medskip
\noindent
\begin{theorem}[Cylinder Uniformization Theorem] \label{T:Unif}
 Let $S$ be a carpet of measure zero in $\Sph^2$ whose peripheral circles are uniformly 
relatively separated uniform quasicircles, and $C_1$ and $C_2$ be distinct peripheral circles of $S$. Then there exists a quasisymmetric map $f$ from $S$ onto a square carpet $T$ in a $\C^*$-cylinder $A$ such that $f(C_1)$ is the inner boundary component of $A$ and $f(C_2)$ is the outer one.
\end{theorem}

This is  \cite[Theorem~1.6]{mB04}. In this statement,  $S$ is  equipped with the spherical metric as by our convention adopted in the introduction.   For the metric on $T$ one can choose the spherical metric, the Euclidean metric, or flat metric on $\C^*$; they are all comparable on $A$ and hence on $T$. 

\medskip
 Let $S$ be a carpet in $\Sph^2$  and $C_1$ and $ C_2$ be two distinct peripheral circles of $S$ that bound the complementary components $D_1$ and $D_2$ of $S$, respectively. We denote by $\Ga(C_1,C_2;S)$ the family of all open paths $\ga$ in $\Sph^2\setminus(\overline D_1\cup\overline D_2)$ that connect $\overline D_1$ and $\overline D_2$.

The following proposition gives an explicit description 
for the extremal mass distribution for the carpet modulus 
$\Mod_S(\Ga(C_1,C_2;S))$ under suitable conditions.

\begin{proposition}\label{P:Expldescr} Let $S$ be a carpet of measure zero in $\Sph^2$ whose peripheral circles are uniformly 
relatively separated uniform quasicircles, and $C_1$ and $C_2$ be two distinct peripheral circles of $S$. 
Then an extremal mass distribution $\rho$ for $\Mod_S(\Ga(C_1,C_2;S))$ exists, has finite and positive total mass, and is given as follows:
Let $f$ be a quasisymmetric map of $S$ to  a square carpet $T$ in a $\C^*$-cylinder $A=\{z\in\C\co r\leq|z|\leq R\}$  such that $C_1$ corresponds to the inner and $C_2$ to the outer boundary components of $A$. Then $\rho(C_1)=\rho(C_2)=0$, and for the peripheral circles $C\ne C_1, C_2$ of $S$ we have 
$$\rho(C)=\frac {\ell_{\C^*}(Q)}{\log(R/r)} ,$$ where $Q$ is the $\C^*$-square bounded by $f(C)$. 
\end{proposition}

This is \cite[Corollary~12.2]{mB04}. Note that a map $f$ as in this proposition exists by the previous Theorem~\ref{T:Unif}.  The map 
$f$ is actually unique up to scaling and rotation around $0$ as follows from Theorem~\ref{P:ACyl}, which  we  will prove in 
Section~\ref{S:Corners}. 
It follows from Proposition~\ref{P:Qcext} that $f$ has a quasiconformal extension to $\Sph^2$, and so $T$ is also a set of measure zero~~\cite[Ch.~IV, \S 1, Section~1.4, Theorem~1.3, p.~165]{LV}. 
From the explicit description of the extremal mass distribution it 
follows that 
$$ 0< \Mod_S(\Ga(C_1,C_2;S))=\frac {2\pi}{\log(R/r)} <\infty.$$

In \cite{mB04} the proof of Proposition~\ref{P:Expldescr} was fairly straightforward, but had to rely on  Theorem~\ref{T:Unif}, whose   proof was rather involved. The only consequence  of Proposition~\ref{P:Expldescr} that we will use is that for the extremal density $\rho$ we have $\rho(C)>0$ for all peripheral circles $C\ne C_1,C_2$. It is an interesting question whether a direct proof of this statement can be given without resorting to
Theorem~\ref{T:Unif}.



\section{Distinguished pairs of peripheral circles}\label{S:Pair}

\no In the following it is convenient to use the term  {\em  square} also for  the boundary of a solid 
Euclidean square  in the usual sense. It will be clear from the context 
which meaning of square is intended.  If $Q$ is a square in either sense, then we denote by $\ell(Q)$ its Euclidean side length.
A {\em corner} of $Q$ is an endpoint  of a side of $Q$.  

 With our terminology we can refer to  the  peripheral circles of the standard Sierpi\'nski carpets $S_p$, $p\ge 3$ odd, simply  as squares. These squares arising as peripheral circles of $S_p$ form a family of uniform quasicircles in the Euclidean metric, since each of them can be mapped to the boundary $\partial Q_0$ of the solid unit square $Q_0$ by a Euclidean similarity. 
 This family is also uniformly relatively separated in the Euclidean metric, because 
 if $C$ and $C'$ are two distinct squares in this family, then for their Euclidean distance we have 
 \begin{eqnarray*}
 \dist(C,C')&\ge& \frac{p-1}2 \min\{\ell(C), \ell(C')\}\\
 &=&\frac{p-1}{2\sqrt2} \min\{\diam(C), \diam(C')\}. 
 \end{eqnarray*}
 Since Euclidean distance and spherical distance on $S_p\sub \Sph^2\cong \widehat \C$ are comparable, it follows that 
 the family of peripheral circles is uniformly relatively separated 
 and consists of uniform quasicircles also with respect to the spherical metric. Since $S_p$ has measure zero in addition, we can apply to $S_p$ the results that were stated in  Section~\ref{S:Aux}.  Moreover, by the comparability of  Euclidean and  spherical metric on $S_p$ the class of quasisymmetric self-maps on $S_p$ is the same for both metrics.

Our goal in  this section is  to prove that any quasisymmetric self-map of $S_p$ preserves the outer and the middle squares as a
pair. By definition the \emph{outer square} $O$ is the peripheral circle that 
corresponds to the boundary of the original unit square in the construction of $S_p$. The \emph{middle square} $M$ is the boundary of the  open middle square  removed from the unit square in the first step of the construction of
$S_p$. It is the unique peripheral circle different from $O$ that is invariant under all Euclidean  isometries of $S_p$. Note that these isometries of 
$S_p$ form  a dihedral group  with eight elements.

\begin{lemma}\label{L:Pair} Let $p\ge 3$ be odd, and  
let $C, C'$ be any (unordered) pair of distinct peripheral circles of $S_p$ other than $M,O$.
Then 
$$
{\rm mod}_{S_p}(\Ga(C,C';S_p))<{\rm mod}_{S_p}(\Ga(M,O;S_p)).
$$
\end{lemma}

\begin{figure}
[htbp]
\begin{center}
\includegraphics[height=40mm]{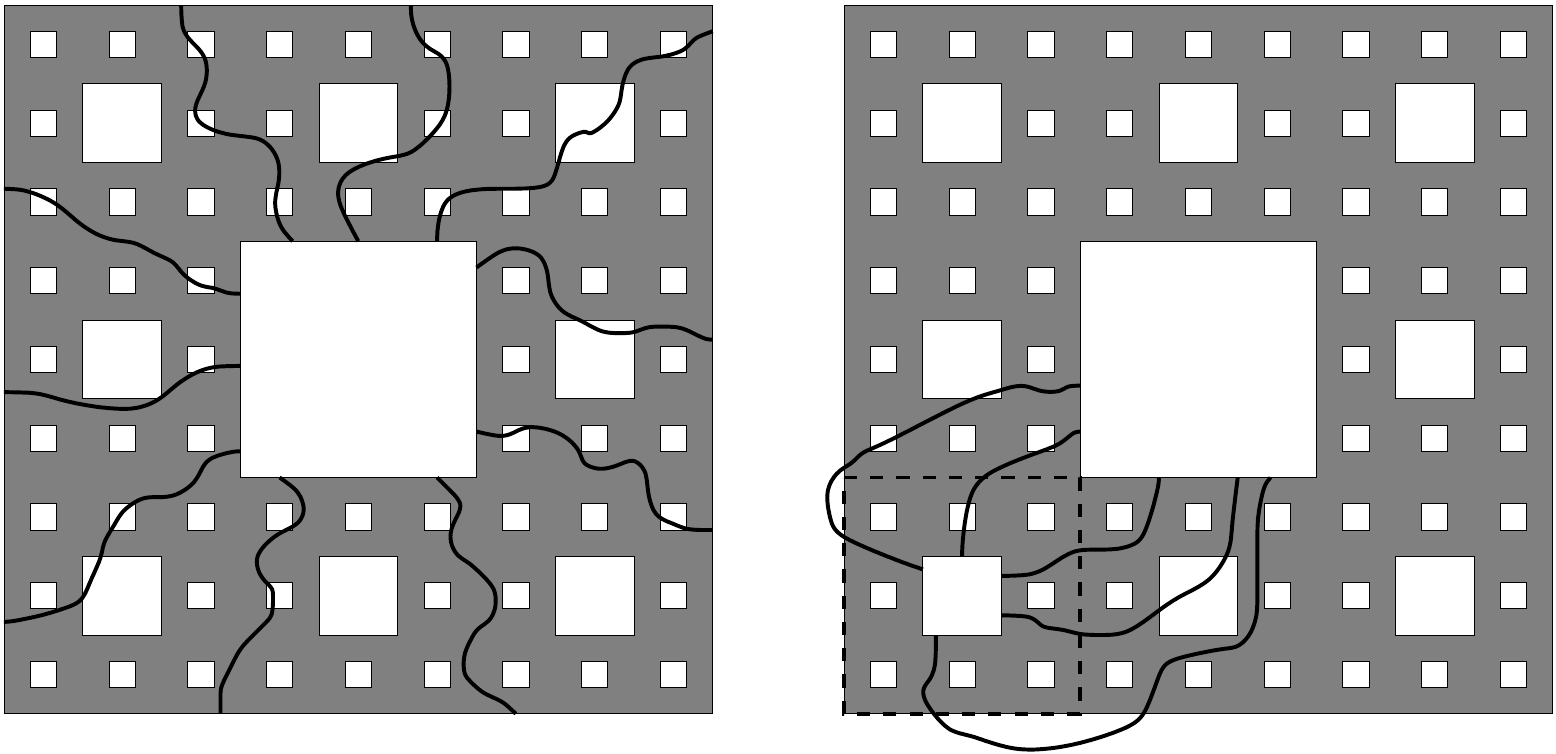}
\caption{
The part of the carpet on the right bounded by the dashed line is a rescaled copy $S$ of $S_p$. 
}
\label{F:Pairofsquares}
\end{center}
\end{figure}

\noindent
\emph{Proof.} The self-similarity of the carpet $S_p$ and the monotonicity property of the modulus give
\begin{equation}\label{eq:nonstrict}
{\rm mod}_{S_p}(\Ga(C,C';S_p))\leq{\rm mod}_{S_p}(\Ga(M,O;S_p)).
\end{equation} 
Indeed, if $l$ and $l'$ are the side lengths of the squares $C$ and $C'$, respectively, then we may assume that $l\le l'$. Then $l\le 1/p^2$. This implies that 
 there exists a copy  $S\subset  S_p$, $S\ne S_p$,  of $S_p$,  rescaled  by the factor $pl$,  so that  $C$ corresponds to $M$, the middle square. Then the  outer square $o$ of $S$ is the rescaled copy of $O$, and the interior region of $o$ is disjoint from $C'$. Hence  every path in $\Ga(C,C';S_p)$ meets $o$ (possibly in one of its endpoints) and so contains a sub-path in $\Ga(C,o;S)$ (see Figure~\ref{F:Pairofsquares} for an illustration of this situation).
Therefore, 
$$
{\rm mod}_{S_p}(\Ga(C,C';S_p))\leq{\rm mod}_{S_p}(\Ga(C,o;S)).
$$
On the other hand, 
$$
{\rm mod}_{S_p}(\Ga(C,o;S))={\rm mod}_{S}(\Ga(C,o;S)),$$
since every  path in $\Ga(C,o;S)$ meets exactly the same peripheral circles of 
$S$ and $S_p$. Moreover, 
$${\rm mod}_{S}(\Ga(C,o;S))={\rm mod}_{S_p}(\Ga(M,O;S_p)),
$$ by Lemma~\ref{L:Invmod}. Inequality  
 \eqref{eq:nonstrict} follows.

To reach a contradiction, assume now that
\begin{equation}\label{eq:modeq} 
{\rm mod}_{S_p}(\Ga(C,C';S_p))={\rm mod}_{S_p}(\Ga(M,O;S_p)).
\end{equation} 
Note that all carpet moduli considered above  are finite by Proposition~\ref{P:Expldescr}, and so  an 
extremal mass distribution exists for each of them by Proposition~\ref{P:infmin}. Then  \eqref{eq:modeq}, the preceding discussion, and the  uniqueness of the extremal mass distributions  implies  that the extremal mass distribution for ${\rm mod}_{S_p}(\Ga(C,C';S_p))$ is obtained from the extremal mass distribution for  
${\rm mod}_{S_p}(\Ga(M,O;S_p))$ by ``transplanting" it to $S$ using a suitable Euclidean similarity between $S$ and $S_p$ (similarly as in the proof of Lemma~\ref{L:Invmod}). Hence the extremal mass distribution for ${\rm mod}_{S_p}(\Ga(C,C';S_p))$ is supported only on the set of peripheral circles of $S_p$ that are also peripheral circles of $S$. 
This is however not the case as follows from 
Proposition~\ref{P:Expldescr}, and we arrrive at a contradiction.  \qed\medskip 


\begin{corollary}\label{C:Group} Let $p\ge 3$ be odd. Then 
every quasisymmetric self-map of  $S_p$ 
preserves the middle  and the outer squares $M$ and $O$ as a pair. 
\end{corollary}

So if $f\co S_p\ra S_p$ is a quasisymmetric map, then $\{f(M),f(O)\}=\{M,O\}$.
This allows the possibility that $f$ interchanges $M$ and $O$, i.e., that 
$f(M)=O$ and $f(O)=M$. We will later see that actually $f(M)=M$ and $f(O)=O$ (Lemma~\ref{L:f(O)=O}). 

\medskip 
\noindent
\emph{Proof.} Assume that $f$ maps the pair $M,O$ to some  pair of peripheral circles $C, C'$.
By Proposition~\ref{P:Qcext}, the map $f$ extends to a quasiconformal homeomorphism  $F$ on $\Sph^2$. 
In particular, $\Ga(C,C';S_p)=F(\Ga(M,O;S_p))$. 
Lemma~\ref{L:Invmod} then  implies 
\begin{equation}\label{E:2Pairs}
{\rm mod}_{S_p}(\Ga(C,C';S_p))={\rm mod}_{S_p}(\Ga(M,O;S_p)).
\end{equation}
By Lemma~\ref{L:Pair} this is only possible if $\{C,C'\}=\{M,O\}$. 
\qed\medskip 

\begin{corollary}\label{C:Gpfinite}
Let $p\ge 3$ be  odd. Then the group ${\rm QS}(S_p)$ 
of quasisymmetric self-maps of  $S_p$ is finite. 
\end{corollary}

\no
\emph{Proof.}
According to Corollary~\ref{C:Group}, the middle  square $M$ and the outer   square $O$ of $S_p$ are preserved as a pair under every quasisymmetric self-map of $S_p$. 
Moreover, by  Corollary~\ref{C:CC}  the group $G$ of all orientation-preserving quasisymmetric self-maps $f$ of $S_p$ 
with $f(M)=M$ and $f(O)=O$ is finite cyclic. If $f_1, f_2\in {\rm QS}(S_p)$ are  orientation-reversing, then $f_1^{-1}\circ f_2\in {\rm QS}(S_p)$ is orientation-preserving. Likewise, if $f_1,f_2\in {\rm QS}(S_p)$  interchange $M$ and $O$, then $ f_1^{-1}\circ f_2$ preserves both $M$ and $O$  setwise. This implies that $G$ is a subgroup of ${\rm QS}(S_p)$ with index at most $4$. Since 
$G$ is a finite group, ${\rm QS}(S_p)$ is  finite  as well.\qed

\section{Quasisymmetric rigidity of square carpets}\label{S:Corners}

\no 
In this section  we prove quasisymmetric rigidity  results for  square carpets  in rectangles and for square carpets  in  $\C^*$-cylinders.  

\begin{figure}
[htbp]
\begin{center}
\includegraphics[height=40mm]{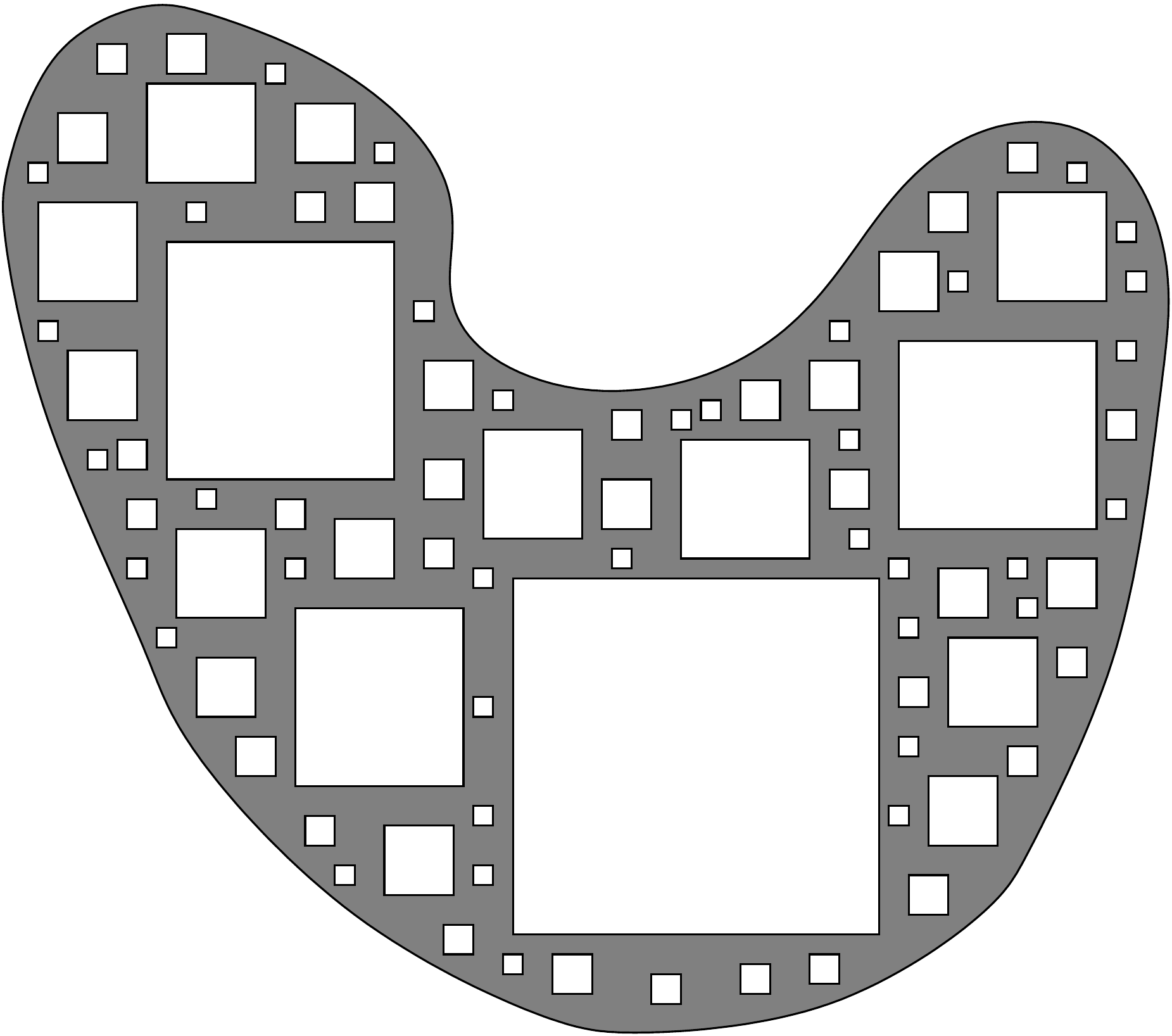}
\caption{
A square carpet in a closed Jordan region.
}
\label{F:Sqcrpt}
\end{center}
\end{figure}

By definition 
 a \emph{square carpet} $S$ {\em in a closed  Jordan region} $D\sub \R^2\cong\C$ is a carpet  $S\sub  D$ so that $\partial D$ is a peripheral circle of $S$ and all other peripheral circles are squares with sides parallel to the coordinate axes (see Figure~\ref{F:Sqcrpt}). We equip such a carpet  with the Euclidean metric.

We will now prove Theorem~\ref{P:Cor}. In the ensuing proofs all metric concepts refer to the Euclidean metric on  $\R^2\cong\C$. 
Moreover, we will use the Euclidean metric also as a base metric in the definition of conformal modulus of a path family. 

\medskip 
\noindent
\emph{Proof of Theorem~\ref{P:Cor}.} Without loss of generality we may assume that $\tilde a\leq a$. Suppose $f$ is a quasisymmetric map as in the statement. Note that the peripheral circles of $S$ are uniform 
quasicircles; this is clear if we use the Euclidean metric, but 
on $K$ the Euclidean and the spherical metrics are comparable, and so the peripheral circles of $S$ are uniform quasicircles also 
with respect to the spherical metric.
Hence  by Proposition~\ref{P:Qcext} the  map $f$ extends to a quasiconformal map $F$ on $\Sph^2$.

 We denote by $C_i$, $i\in \N$,  the peripheral circles of $S$  distinct from $\partial K$, and by $Q_i$ the closed solid square bounded by $C_i$.   We set $\widetilde C_i:=f(C_i)$. Then the sets $\widetilde C_i$, $i\in \N$, 
are the  peripheral circles of $\widetilde S$  distinct from $\partial \widetilde  K$. For $i\in \N$  let $\widetilde Q_i$ be the closed solid square 
bounded by  $\widetilde C_i$. Note that then $\widetilde Q_i=F(Q_i)$. 

For $t\in [0,a]$ we  
denote by $\ga_t$ the path $u\mapsto{t}+\iu u$, $u\in[0,1]$, and let $\Gamma=\{\ga_t:t\in [0,a]\}$. So $\Gamma$ consists precisely  of the closed  vertical  line
segments connecting the horizontal sides of $K$.
Then $\Mod(\Gamma)=a$ and  the function $\rho_0\equiv 1$ on $K$ is an extremal 
density for $\Mod(\Gamma)$. If  $\tilde \rho$ is another extremal density for $\Mod(\Gamma)$, then $\tilde \rho(z)=1$ for almost every $z\in K$.  

We define a Borel density $\rho_1$ on $K$ as follows. For $z\in K$ we set 
 $\rho_1(z)=\ell( \widetilde Q_i)/\ell(Q_i)$ if $z\in Q_i$ for some 
 $i\in I$, and $\rho_1(z)=0$ otherwise.  We will   show that $\int_\gamma\rho_1\,ds\ge 1$ for almost every path 
 $\gamma\in \Gamma$. This will make it possible to adjust $\rho_1$  on a set of measure zero  so that  the resulting density 
 is admissible for $\Gamma$.   

To see this,   let 
 $E$ be the set of all $t\in [0,a]$ 
 for which  $\ga_t\cap S$ has positive length, i.e., 
 $E=\{ t\in [0,a] : \int_{\ga_t}\chi_S\, ds>0\}. $ 
Since $S$ has measure zero, it follows from  Fubini's theorem that $E$ has   $1$-dimensional  measure zero. 

Moreover, since quasiconformal maps are absolutely 
continuous on almost every 
line~\cite[Theorem~3.1, p.~170]{LV},
there exists a set 
$E'\sub [0,a]$ of $1$-dimensional measure zero such that the map $u\mapsto F(t+\iu u)$ is  absolutely continuous on $[0,1]$ for each   $t\in [0,a]\setminus E'$. 

If we define $E_0=E\cup E'\sub [0,a]$, then $E_0$ also has $1$-dimensional measure zero.
Moreover, if $t\in [0,a]\setminus E$, then the  map $u\mapsto F(t+\iu u)$ is  absolutely continuous on $[0,1]$, and $\ga_t\cap S$ has length zero. It follows that 
 $F(\ga_t\cap S)= F(\ga_t)\cap 
\widetilde S$  also has  length zero. By our normalization assumption, the map  $f$, and hence also its extension $F$,   sends each  horizontal side of $K$ to a horizontal side of $\widetilde K$. Thus $F(\ga_t)$ joins the horizontal sides of $\widetilde 
K$ and we conclude that 
$$\int_{\gamma_t}\rho_1\, ds = \
\sum_{\ga_t\cap Q_i\ne \emptyset}\ell(\widetilde Q_i)=
\sum_{F(\ga_t)\cap\widetilde Q_i\ne \emptyset}\ell(\widetilde Q_i)\ge 1. $$

For $z\in K$ define  $\rho_2(z)=\infty$ if $z\in E_0\times [0,1]$ and 
$\rho_2(z)=0$ otherwise. Then 
$\int_{\gamma_t}\rho_2\, ds =\infty$ for $t\in E_0$. It follows that   $\rho=\rho_1+\rho_2$ is admissible for $\Gamma$. 
Moreover, since $\rho_2(z)=0$ for almost every $z\in K$, we have  
$$\int_K\rho^2\, dA=\int_K\rho_1^2\, dA=\sum_i\ell(\widetilde Q_i)^2=\widetilde a\le a=\Mod(\Gamma),$$
 where $dA$ indicates integration with respect to $2$-dimensional Lebes\-gue measure. 
Therefore,  $\widetilde a=a$, $\widetilde K=K$,  and  $\rho$ is extremal for 
$\Mod(\Gamma)$. Hence $\rho(z)=1$ for almost every $z\in K$ which in turn implies that $\ell(Q_i)=\ell(\widetilde Q_i)$ for all $i\in \N$. So  each square $Q_i$ has the same side length as its image square  
$\widetilde Q_i=F(Q_i)$.

We will next show that actually  $Q_i=\widetilde Q_i$ for each $i\in \N$. To see this, let $i\in \N$ be arbitrary, and  
consider  the family $\Gamma'$ of 
all  open   line segments   that are  parallel to the real axis and connect  the left vertical  side of $K$ to the right vertical  side of $Q_i$.   Let  $\ga\in \Gamma'$ and assume  in addition that 
$\ga\cap S$ has length zero and  that 
$F$ is absolutely continuous on $\ga$. Then  the intersection $F(\ga)\cap \widetilde S$
also has  length zero, and we have 
$${\rm length} (\ga)=\sum_{\ga\cap Q_j\ne \emptyset}\ell(Q_j).  $$ Since $F$ sends each peripheral square
  to a 
square of the same side length, we conclude that for the length $L$ of
the projection of $F(\ga)$ to the real  axis 
we have 
$$ L\le \sum_{F(\ga)\cap \widetilde Q_j\ne \emptyset} \ell (\widetilde Q_j)=\sum_{\ga\cap Q_j\ne \emptyset}\ell(Q_j)={\rm length} (\ga). $$
In other words, the length of the projection of $\ga$ to the real axis (which is equal to the length of $\ga$) does not increase under the map $F$.

By an argument similarly as above, one can see that the family 
$\Gamma_0'$ of all line segments $\ga\in \Gamma'$ for which $\ga\cap S$ has positive length or for which $F$ is not absolutely continuous on $\ga$ has modulus zero. 
In particular, for each  $\ga\in \Gamma'$ there are   line 
segments $\ga'\in \Gamma'\setminus \Gamma_0'$ that are arbitrarily close to $\ga$. Since the length of the projection of each such line segment $\ga'$ to the real axis does not increase under $F$, it  cannot increase for any $\ga\in \Gamma'$  either.

 We conclude that the distance of $\widetilde
Q_i$ to the left vertical  side of $\widetilde K=K$ is bounded from above by the distance of $Q_i$ to the 
left vertical  side of $K$.
Using the same argument for the inverse map, we
conclude that these distances are actually equal.  
We can apply the same reasoning  for   other pairs of respective sides  of $K$ and $Q_i$. 
Hence  $Q_i$ and $\widetilde Q_i$ are squares in $K=\widetilde K$ with the same distances to all sides of $K$. This implies     $Q_i=\widetilde Q_i$.

Moreover, we can actually deduce  that $F$ maps each side of  $C_i=\partial Q_i$ into itself. Indeed, let $p\in C_i$ be a point on the right vertical  side of $C_i$, say. If $\Gamma'$ is as above, then there is a line segment $\ga\in \Gamma'$ which has $p$ as one of its endpoints. 
Then $F(p)\in \widetilde C_i=C_i$ is one endpoint of $F\circ \ga$, while the other endpoint lies on the left vertical  side of 
$K$. As we have seen, the length of the projection of $F\circ \ga$ to the real axis is bounded by the length of $\ga$ which is equal to the distance of $C_i$ to the left vertical  side of $K$. This is only possible if $F(p)$ lies on the left vertical side of $C_i$. So $F$ maps this side into itself, and,  similarly,   each  side  of $C_i$ into itself. 
This in turn implies $F$ must fix each  corner  of $C_i$.  

Since $i\in \N$ was arbitrary,  we conclude $S=\widetilde S$, and that $F$, and hence also $f$, fixes the corners  of all squares $Q_i$, $i\in \N$. 

Now every  subset of  a carpet that meets all but finitely many peripheral circles is dense. In particular, the set $D$ consisting of all corners  of the squares 
$C_i$, $i\in \N$, is a dense subset of $S$. Since $f$ is the identity  on $D$, it follows that 
 the map $f$ is the identity on $S$. \qed
\medskip 

We now prove Theorem~\ref{P:ACyl}. The argument  is very similar to the proof of 
Theorem~\ref{P:Cor}. In the proof metric notions refer to  the  flat metric on $\C^*$ given by the length element $|dz|/|z|$. We will also use  it  as a  base metric for modulus.  For terminology related to  $\C^*$-cylinders and $\C^*$-squares used in the ensuing proof see the discussion before Theorem~\ref{T:Unif}.

\medskip 
\noindent
\emph{Proof of Theorem~\ref{P:ACyl}.} Let  $f$ be  as in the statement.  Each $\C^*$-square $Q$ (equipped with the flat metric 
on $\C^*$) that satisfies  $\ell(Q)\le \pi$ is isometric to a Euclidean square of the same sidelength. 
So the family of all peripheral circles of $S$ that bound such $\C^*$-squares consists of uniform quasicircles. There are only finitely many peripheral circles not in this family, namely the boundary components of $A$ and   boundaries of complementary components of $S$ that are $\C^*$-squares $Q$ with $\pi<\ell(Q)<2\pi$. Since each of these finitely peripheral circles  (equipped with the flat metric
 on $\C^*$) is bi-Lipschitz equivalent to the unit circle, it follows that the family of {\em all}  peripheral circles of $S$ consists of uniform quasicircles. So again by  Proposition~\ref{P:Qcext} we can extend   the map $f$ to a quasiconformal map $F$ on $\Sph^2$. Then $F$, as the map $f$,  is   orientation-preserving.

We may assume that the inner components of $A$ and $\widetilde A$ are equal to the unit circle, and the outer boundary components of $A$ and $\widetilde A$ are equal to $\{z\in \C:|z|=R\}$ and $\{z\in \C:|z|=\widetilde R\}$, respectively, where  $1<R\le  \widetilde R$.  We denote by $Q_i$, $i\in \N$,  the (open) $\C^*$-squares
whose boundaries give  the peripheral circles  of $S$  distinct from the boundary components of $A$, and set $\widetilde Q_i=F(Q_i)$ 
for $i\in \N$. Then $\partial \widetilde Q_i$, $i\in \N$, is the family of all peripheral circles of $\widetilde S$ distinct from the boundary components of $\widetilde A$.

We now consider the family $\Gamma$ of closed radial line segments
joining  the boundary components of $A$. Then 
$\Mod(\Gamma)= 2\pi/\log(R)$, and $\rho_0=1/\log(R)$ is the essentially unique extremal density (with the flat metric as the underlying base metric).
On the other hand, we define a density $\rho_1$ on $A$ such that   
$$\rho_1(z)=\frac{\ell_{\C^*}(\widetilde Q_i)}{\log(\widetilde R)\ell_{\C^*}(Q_i)}$$
if 
 $z\in Q_i$ for some $i\in \N$, and $\rho(z)=0$ elsewhere on $A$. As in 
proof of 
Theorem~\ref{P:Cor}, one shows that up to adjustment on a set of measure zero, $\rho_1$ is admissible for $\Gamma$. 
Moreover, 
$$\int_A \rho_1^2\, dA_{\C^*}=\frac 1{\log^2(\widetilde R)} \sum_i  
\ell_{\C^*}(\widetilde Q_i)^2=\frac {2\pi }{\log(\widetilde R)} \le 
\frac {2\pi }{\log( R)}=\Mod(\Gamma),$$
where $dA_{\C^*}$ means integration with respect to the area element induced by the flat metric. 
Hence $R=\widetilde R$, $A=\widetilde A$, and  $\rho_1$ (up to a change on a set of measure zero) is extremal for $\Mod(\Gamma)$. We  conclude that 
$\rho_1=1/\log(R)$ almost everywhere on $A$, which implies that  $\ell_{\C^*}(Q_i)=\ell_{\C^*}(\widetilde Q_i)$ for all $i\in \N$. 
So again each $\C^*$-square $Q_i$ has the same side length as its image $\widetilde Q_i$ under $F$.

Using this and 
arguments similar to the ones in the proof of Theorem~\ref{P:Cor},  one can show that for each $i\in \N$ the  squares  $Q_i$ and $\widetilde Q_i$ have the same 
distances  to  the inner and outer boundary components of $A$, and that $F$ maps the bottom and top sides 
of $Q_i$ to the bottom and top sides of $\widetilde Q_i$, respectively. This implies that $F$ sends the corners of $Q_i$ to the corners of $\widetilde Q_i$. Since  $F$  is orientation-preserving, the cyclic order of the corners is preserved under the map $F$.  It follows that 
for each $i\in \N$, there exits a rotation $r_i$ around $0$ such that $r_i(Q_i)=F(Q_i)=\widetilde Q_i$, and such that $r_i(c)=f(c)=F(c)$ if $c$ is a corner of $Q_i$. 

 So far, we exclusively used the behavior of $F$ on radial directions. We will now investigate the behavior of $F$ on  ``circular directions". To do this, we consider the circular projection 
$$P_i:= \{t\in \R: te^{\iu\alpha}\in Q_i \text{ for some } \alpha\in [0,2\pi]\}$$  of $Q_i$ to the real axis.
Each set $P_i$, $i\in \N$,  is an open subinterval of $(1,R)$. 

If 
 $P_i\cap P_j\ne \emptyset$ for $i,j\in \N$, $i\ne j$, then the circular projections of the    squares  $Q_i$ and $Q_i$ overlap, and so we can find a family $\Gamma'$ of closed  circular arcs $\ga$,  each contained in a circle of radius  $t\in P_i\cap P_j$, that  join  two vertical sides of $Q_i$ and $Q_j$ facing each other. This family $\Gamma'$ has positive modulus, and similarly as in the proof of Theorem~\ref{P:Cor},  one can show that the length of the radial projection of each path $\gamma\in \Gamma'$   to the unit circle does not increase under the map $F$.  Applying the same argument to the other pair of   vertical sides of $Q_i$ and $Q_j$ facing each other, we conclude that   the circular distance of $Q_i$ and $Q_j$ is the same as the circular distance of image squares $\widetilde Q_i$ and $\widetilde Q_j$. This implies that  $r_i=r_j$.

We can write 
$$U:=\bigcup_iP_i=\bigcup_{k\in J}M_k,$$
where $J$ is a countable index set, and the sets $M_k$, $k\in J$,  are pairwise disjoint open subintervals of $(1,R)$ forming the connected components of $U$. Suppose that $M_k=(a_k,b_k)$, where $1\le a_k<b_k\le R$, and set 
$A_k:=\{ z\in \C: a_k<|z|<b_k\}$.
For every $i\in \N$ there exists precisely one $k\in J$ such that 
$Q_i\sub A_k$. Moreover, since any two points $u,v\in M_k$ can be connected by a chain of intervals $P_i\sub A_k$, it follows that 
$r_i=r_j$ whenever $P_i,P_j\sub M_k$ for some $k\in J$.
So  for each $k\in J$ there exists a rotation $\tilde r_k$ around $0$ such that $r_i=\tilde r_k$ whenever $Q_i\sub A_k$.

We claim that 
\begin{equation} \label{eq:subann}
f|S\cap A_k=\tilde r_k|S\cap A_k
\end{equation} 
for each $k\in J$. To see this let $k\in J$ and $z_0\in S\cap A_k$ be arbitrary. Since the set of corners of the $\C^*$-squares $Q_i$ is dense in $S$, there exists a sequence $(c_n)$ of such corners  
with $c_n\to z_0$ as $n\to \infty$. If $c_n$ is a corner of the $\C^*$-square $Q_{i_n}$, then $Q_{i_n}\sub A_k$ for sufficiently 
large $n$. For these $n$ we have $\tilde r_k(c_n)=r_{i_n}(c_n)=f(c_n)$.
Passing to the limit $n\to \infty$, we conclude that indeed 
$\tilde r_k(z_0)=f(z_0)$ as desired. 

We also have  
\begin{equation}\label{eq:rotonQ} 
 f|\partial Q_i=r_i|\partial  Q_i
\end{equation}
for each $i\in \N$. Indeed, let $i\in \N$ be arbitrary. 
By \eqref{eq:subann} it is clear that $f$ and $r_i$ agree 
on the interior of each vertical side of $Q_i$, because these interiors are contained in a suitable set $S\cap A_k$. 
Let $u$ be a point on one of the other sides of $Q_i$, say on the bottom side of $Q_i$. Pick a corner $v$ of $Q_i$ on the same   side. We will construct a sequence $(k_j)$ in $J$, and sequences 
$(u_j)$ and $(v_j)$ of points such that $u_j,v_j\in S\cap  A_{k_j}$ for $j\in \N$, and $u_j\to u$ and $v_j\to v$ as $j\to \infty$. 
Then by  \eqref{eq:subann}, we have $\tilde r_{k_j}(u_j)=f(u_j)$ and 
$\tilde r_{k_j}(v_j)=f(v_j)$. By passing to a subsequence if necessary, we may assume that the rotations $r_{k_j}$ converge to a rotation $r'$ uniformly on $A$ as $j\to \infty$. 
Then $r'(u)=f(u)$ and $r'(v)=f(v)$. Since $v$ is a corner of $Q_i$, we also have $r_i(v)=f(v)=r'(v)$. Hence $r'=r_i$, and so $r_i(u)=r'(u)=f(u)$ as desired. 

To produce the sequences $(u_j)$ and $(v_j)$, we  consider the set of 
$E=[1,R]\setminus U$. This is the set of all radii of circles centered at $0$ that lie in $S$. Since $S$ has measure zero, 
$E$ has $1$-dimensional measure zero.

Suppose that $u=se^{\iu\alpha}$,
where $1<s=|u|=|v|<R$ and $\alpha\in [0,2\pi]$. Since $E$ is a set of measure zero, we can find a sequence $(s_j)$ of ``good radii" 
such that $s_j\in (1,s)\setminus E$ and $s_j\to s$ as 
$j\to \infty$.  Then there exists $k_j\in J$ such that 
$s_j\in A_{k_j}$ for  $j\in \N$. 
Define  $u'_j=s_je^{\iu\alpha}$ for $j\in \N$. Then 
$u'_j\in A_{k_j}\sub A$ for $j\in \N$, and  $u'_j\to u$ as $j\to \infty$, but  
the sequence $(u'_j)$ is not necessarily contained in $S$. To achieve this we shift each point $u'_j$ on the circle $\{z: |z|=s_j\}$
if necessary. More precisely, 
if $u'_j\in S$, we let $u_j:=u'_j$.
If $u'_j$ does not lie in $S$, then $u'_j$ is contained in one of the $\C^*$-squares $Q_l$, $l\in \N$. We can then move $u'_j$ on the circle $\{z: |z|=s_j\}$ to a point $u_j$ on one of the  vertical  sides of $Q_l$. Note that $Q_l\ne Q_i$, and so the diameter of $Q_l$ is small if $j$ is large, since $\C^*$-squares  $Q_l\ne Q_i$ exceeding a given size cannot be arbitrarily close to $Q_i$. 
So we have  $u_j\in S\cap A_{k_j}$ and $u_j\to u$ as $j\to \infty$.
A  sequence $(v_j)$ with $v_j\in S\cap  A_{k_j}$ for $j\in \N$ and with $v_j\to v$ as $j\to \infty$ is constructed similarly. Note that our construction guarantees that  $u_j$ and $v_j$ lie in $S$ and in  the same set  $A_{k_j}$ which was  crucial for the argument   
in the previous paragraph. 

Now that we have established \eqref{eq:rotonQ},  we can finish the argument as follows. The proof of Proposition~5.1 in \cite{mB04} 
combined with  \eqref{eq:rotonQ} shows that one can find a  quasiconformal  extension $\widetilde F$ of $f$ to $\Sph^2$ such that 
$\widetilde F|Q_i=r_i|Q_i$ for each $i\in \N$. Then $\widetilde F$ is conformal on each  $Q_i$. Since  $\widetilde F$ is quasiconformal and the squares   
$Q_i$ fill $A$ up to a set of measure zero, it follows that $\widetilde F$ is a $1$-quasiconformal map on the interior of $A$. Hence $\widetilde F$ is a conformal map  on the interior of 
$A$~\cite[Theorem~5.1, p.~28]{LV}. 
Since $\widetilde F=r_1$ on $Q_1$, it follows that 
$\widetilde F|A$ is a rotation around $0$. Then on $S$ the map
$f=\widetilde F|A$ also agrees with such a  rotation.  The 
statement follows. 
\qed\medskip

\section{Weak tangent spaces}\label{S:Wtangents}\no 
In this section   we discuss some facts about weak tangents of the carpets $S_p$. The most important result here is Proposition~\ref{P:Corwtangents} which will be  crucial  in the proofs of our main theorems. 

In general, weak tangent spaces can be defined as  Gromov-Hausdorff limits of pointed metric spaces obtained by rescaling the underlying metric (see \cite[Chapters 7,8]{BB01} and  \cite[Chapter 8]{DS97} for the general definitions, and \cite{BK02a} for applications very similar in spirit to the present paper).  As we will need this only for the carpets $S_p$, we will first present a suitable definition for arbitrary subsets of $\widehat\C
=\C\cup\{\infty\}\cong\Sph^2$ and then  further adjust the definition for the carpets $S_p$. 

If $a,b\in \C$, $a\ne 0$, and  $M\sub  \widehat\C$, we denote by 
$ a M+b$ the image of $M$ under the M\"obius transformation $z\mapsto a z+ b$ on $\widehat\C$. 
Let $A$ be a subset of $\widehat\C$ with a distinguished point $z_0\in A$, $z_0\ne \infty$. We say that a closed set $A_\infty\subseteq \widehat\C$ is a \emph{weak tangent}
of $A$ {\em (at $z_0$)} if there exists a  sequence  $(\lambda_n)$ 
of positive real numbers with $\lambda_n\to\infty$ such that the  the sets 
$A_n:=\lambda_n(A-z_0)$ tend to $A_\infty$  as $n\to \infty$ in the sense of Hausdorff convergence on $\widehat\C$ equipped with the spherical metric (see  \cite[Chapter~7, Section 7.3.1]{BB01} for the definition of Hausdorff convergence of sets in a metric space). 
In this case, we use the notation
$$
A_\infty=\lim_{n\to \infty}  (A, z_0,\lambda_n).
$$
So  a weak tangent of $A$ at $z_0$ is obtained by 
extracting a limit from  
``blowing up"  $A$ at $z_0$ by suitable scaling factors 
$\lambda_n\to \infty$. 
Every weak tangent of $A$ contains $0$, and, if $A$ is not a singelton set, also the point $\infty$. 

A set $A\sub \widehat\C$ has  weak tangents at each point $z_0\in A\setminus \{\infty\}$, because for every sequence   $(\lambda_n)$ of positive numbers with $\lambda_n\to\infty$, there is a subsequence $(\lambda_{n_k})$ such that the sequence of the sets  $A_{n_k}=\lambda_{n_k}(A-z_0)$ converges as $k\to \infty$
(\cite[Theorem~7.3.8, p.~253]{BB01}).
In general,  weak tangents 
at a point  are not unique. In particular, if $\lambda>0$ and  $A_\infty$ is a weak tangent of $A$ at a point, then $\lambda A_\infty$ is also  a weak tangent.  

It is  advantageous to avoid  this scaling ambiguity 
of weak tangents for  the standard carpets $S_p$, $p\ge 3$ odd,
and restrict the scaling factors  $\lambda_n$  used in the definition of a weak tangent to powers of $p$. So in  the following, {\em a weak tangent of $S_p$  
at a point $z_0\in S_p$} is a closed set $A_\infty\sub \widehat \C$ such that 
$$A_\infty=\lim_{n\to \infty}  (S_p, z_0,p^{k_n}), $$
where $k_n\in \N_0$ and $k_n\to \infty$ as $n\to \infty$. 
If this  limit exists along the full sequence $(p^n)$, i.e., if  
$$A_\infty=\lim_{n\to \infty}  (S_p, z_0,p^{n}) $$ exists, then 
$A_\infty$ is the unique weak tangent of $S_p$ at $z_0$.
We equip each weak tangent of $S_p$ with the spherical metric unless otherwise indicated.

We now exhibit some points in $S_p$, where we have unique weak tangents, and set up some notation for the 
weak tangents thus obtained.  

 \begin{figure}
[htbp]
\begin{center}
\includegraphics[height=60mm]{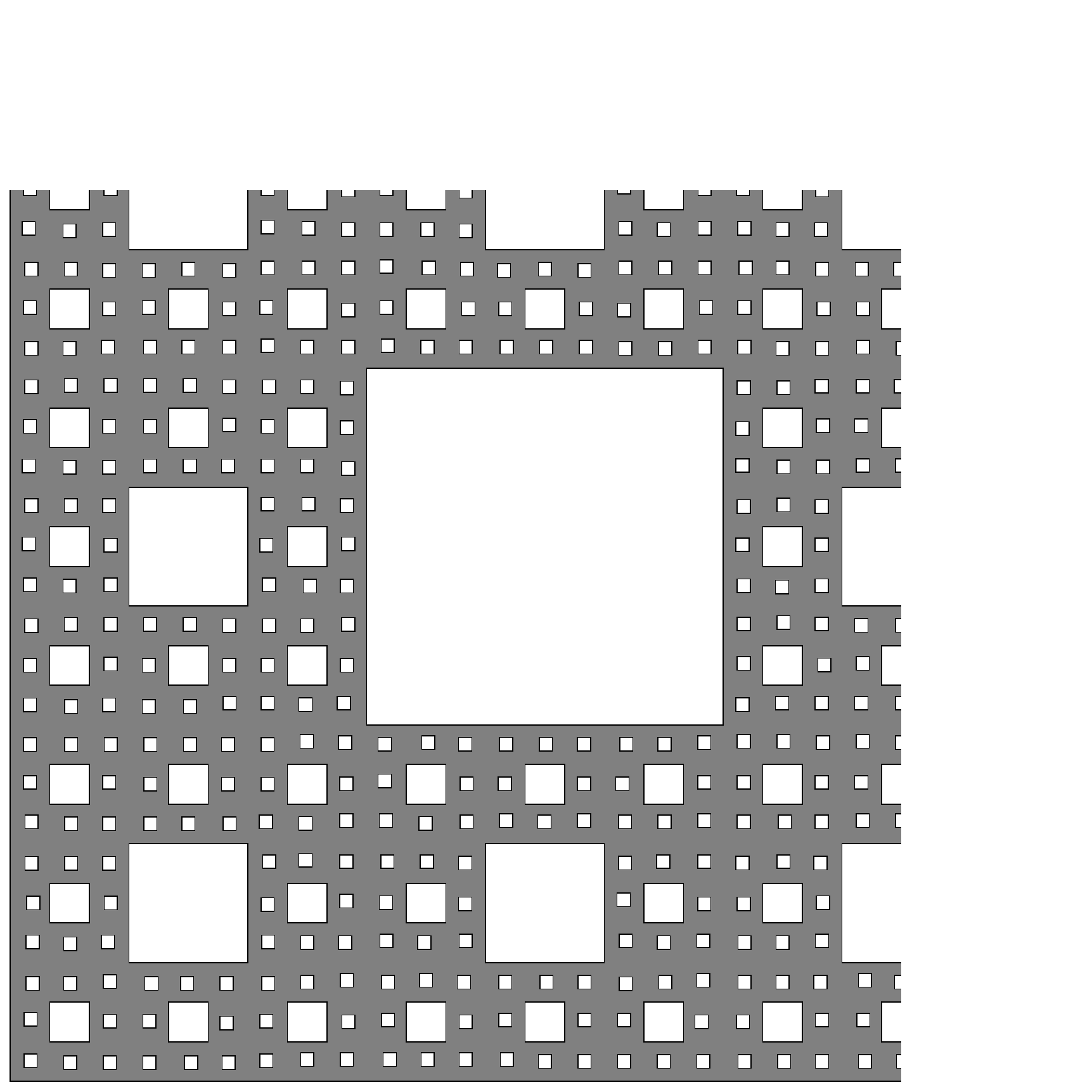}
\caption{
The weak tangent space $W_{\pi/2}$ for $p=3$.
}
\label{F:WT}
\end{center}
\end{figure}

Fix an odd integer $p\ge 3$. At the point $0$ the carpet $S_p$ has the unique weak tangent
\begin{equation}\label{eq:wp2}
W_{{\pi}/2}:=\lim_{n\to \infty}(S_p, 0,p^n)=\{\infty\}\cup \bigcup_{n\in \N_0} p^nS_p.
\end{equation}
The existence of this  and similar limits below can easily be justified  by observing that by self-similarity of $S_p$ the relevant sets involved form an increasing sequence. Here  \eqref{eq:wp2} follows from the inclusions 
$p^nS_p\sub p^{n+1}S_p$ for $n\in\N_0$. 

Similarly,  at each corner of $O$ there exists a unique weak tangent of $S_p$ 
obtained by a  suitable rotation of the set $W_{{\pi}/2}$ around $0$.

Let $m=1/2$ be the midpoint of the bottom side of $O$. Then 
at the point $m$ the carpet $S_p$ has the  unique weak tangent
$$
W_{\pi}:=\lim_{n\to \infty} (S_p, 1/2, p^n)=\{\infty\}\cup \bigcup_{n\in \N_0} p^n(S_p-m).
$$
Moreover, if $z_0$ is any midpoint of a side of a square that is a peripheral circle of $S_p$, then $S_p$ has a unique weak tangent at $z_0$ obtained by a  suitable rotation of the set $W_{{\pi}}$ around $0$. 
This easily follows from  the fact   that 
for the existence and uniqueness of a weak tangent at a point $z_0$, only an arbitrary small (not necessarily open) relative neighborhood of $z_0$ in $S_p$ is relevant, as the complement of such a neighborhood will ``disappear to infinity" if the set is blown up at 
$z_0$.   Moreover, by   self-similarity of $S_p$ for each of these midpoints $z_0$  we can choose a relative neighborhood $N$ of $z_0$ in $S_p$ so that a suitable 
Euclidean similarity maps $N$ to $S_p$ and $z_0$ to $m$, and where the scaling factor of the similarity is an integer power of $p$.

Let 
$$ c:=\frac {p-1}{2p}+ \frac {p-1}{2p} \iu$$
be the lower left corner of $M$. Then at $c$ the carpet 
$S_p$ has a unique weak tangent
$$
W_{{3\pi}/2}:=\lim_{n\to \infty} (S_p,c , p^n)= 
\{\infty\}\cup \bigcup_{n\in \N_0} p^n\big(\iu S_p\cup (-\iu)S_p\cup (-1)S_p\big).
$$
Note that $W_{{3\pi}/2}$ can  be obtained by pasting together three copies
of $W_{{\pi}/2}$.  If $z_0$ is any corner of a square $C\ne O$ that is a peripheral circle of $S_p$, then $S_p$ has a unique weak tangent at $z_0$ obtained by a  suitable rotation of the set $W_{3{\pi}/2}$ around $0$.

The angles $\pi/2$, $\pi$, $3\pi/2$ used as  indices of the above weak tangents  indicate that the corresponding space  is contained in the closure of the quarter, half-, and three-quarter plane, respectively. Mostly, it will be  clear from the context what $p$ is,   
so we will usually omit an additional label $p$ from  the notation; if we want to indicate $p$, then we write
$W_{\pi/2}(p)$ for the weak tangent $W_{\pi/2}$ of $S_p$, etc. 

For $p\ge 3$  odd, we denote by $D_p$ the set of all  midpoints of sides and all corners of squares  that are peripheral circles of $S_p$.
It follows from our previous discussion that at each point $z_0\in D_p$ the carpet $S_p$ has a unique weak tangent $W$ isometric to one of the sets 
$W_{{\pi}/2}$, $W_{\pi}$, or  $W_{{3\pi}/2}$. 
The weak tangent $W$ of $S_p$ at $z_0\in D_p$    can always be written as 
\begin{equation}\label{eq:wunion}
 W=\{\infty\}\cup \bigcup_{n\in \N_0} p^n(N-z_0),
 \end{equation} 
where $N$ is a suitable relative neighborhood of $z_0$ in $S_p$
such that  
\begin{equation}\label{eq:asc}
p^n(N-z_0)\sub  p^{n+1}(N-z_0)
\end{equation}
 for all $n\in \N_0$. 
Actually, we can choose $N$ to be a rescaled copy of $S_p$ 
(if $W$ is isometric to $W_{{\pi}/2}$ or  $W_{\pi}$) or a union of three rescaled copies of $S_p$ (if $W$ is isometric to $W_{{3\pi}/2}$). Note that \eqref{eq:wunion} and \eqref{eq:asc}
 imply that $p^nW =W$ for all $n\in \Z$.

\begin{lemma} \label{L:Wunifqcirc} Let $p\ge 3$ be odd, $z_0\in D_p$, and $W$ be the weak tangent of $S_p$ at $z_0$. 
Then $W$ is a  carpet of measure zero. If $W$ is   equipped with the spherical metric, then 
  the  peripheral circles of $W$ form a family of uniform quasicircles that are uniformly relatively  separated. 
\end{lemma}

\no{\em Proof.} We know that up to rotation around $0$, the set $W$ is equal to one of the weak tangents $W_{{\pi}/2}$, $W_{\pi}$,  $W_{{3\pi}/2}$. So it is enough to show the statement  for these weak tangents. We will do this 
  for $W_{{\pi}/2}$. The proofs for $W_{\pi}$ and   $W_{{3\pi}/2}$ are  the same with minor modifications. 

First note that $W_{{\pi}/2}$ is a carpet, since it can be represented  as in
\eqref{eq:carprep}. Moreover,  this set  has measure zero, because  by \eqref{eq:wunion} it  can be written as a countable union of sets of measure zero. 

Let $$\Omega =\{z\in \C: \Re(z)>0 \   \text{and}\  \Im(z)> 0\}$$ be the open quarter-plane whose closure (in $\widehat \C$) contains  
$W_{{\pi}/2}$. Then $\partial \Omega$ is a peripheral circle of $W_{{\pi}/2}$. Since $\partial\Omega$  can be mapped to the unit circle by a bi-Lipschitz map, this peripheral circle is a quasicircle. 

All other peripheral circles of $W_{{\pi}/2}$ are squares; actually, they are precisely the squares of the form $C'=p^nC$, where $p\in \N_0$ and $C$ is a peripheral circle of $S_p$ different from the outer square $O$.  
As all of these peripheral circles are similar to $O$, and $O$ is bi-Lipschitz equivalent to the unit circle, the peripheral circles $C'\ne \partial \Omega$ of $W_{{\pi}/2}$ are uniform quasicircles in the Euclidean metric.  This is equivalent to  a uniform lower bound  
for certain (metric) cross-ratios of points on these peripheral circles (see~\cite[Proposition~4.4 (iv)]{mB04}). Since cross-ratios for points in $\C$ are the same in the Euclidean and in the chordal metric (the restriction of the Euclidean metric on $\R^3$ to $\Sph^2\cong \widehat \C$), it follows 
 that the peripheral circles $C\ne \partial \Omega$ of $W$ form  a family of uniform quasicircles in the chordal metric. Since  
  chordal and spherical metric on $\widehat \C$ are comparable, we also get a family of  uniform quasicircles in the spherical  metric. If we add  the quasi-circle $\partial \Omega$ to 
 this collection, we still have a family of uniform quasicircles in the spherical metric. 
 
 The uniform relative separation property of the peripheral circles of $W_{{\pi}/2}$ can  be established 
 similarly by passing from the Euclidean to the chordal and the spherical metrics. 
 Namely, first note that if 
 $C$ and $C'$ are peripheral circles of  $W_{{\pi}/2}$ and $C\ne C'$, then 
 for their Euclidean distance we have 
 $$ {\rm dist}(C,C')\ge\frac{p-1}2  \min\{ \ell(C), \ell(C')\}, $$ where 
 $\ell (C)$ and $\ell(C')$ denote the Euclidean side lengths of $C$ and $C'$, respectively, with the convention  $\ell(\partial \Omega)=\infty$. This follows from the fact that if $ \ell(C)\le  \ell(C')$ say, then there exists a rescaled copy $S$ of $S_p$ in $W_{{\pi}/2}$
 with $C\sub S$ such that $C$ corresponds to the middle square 
 of $S_p$, and $C'$ meets $S$ at most in points of the peripheral circle $o$ of $S$ that corresponds to $O$. 
  
The relative uniform separation of the family of all peripheral circles of $W_{{\pi}/2}$ with respect to the Euclidean metric follows. Again this is  equivalent  to  a uniform lower bound for certain metric cross-ratios  (this follows  from 
 \cite[Lemma~4.6]{mB04}; to include $\partial \Omega$, we need a slightly extended form of this lemma where one of the sets is allowed to have infinite diameter, but the statement and the  proof of the lemma can easily be adjusted). Since  cross-ratios are unchanged if we pass to the chordal metric,  it follows that the family of peripheral circles of $W_{{\pi}/2}$ is uniformly relatively separated with respect to the chordal metric, and hence also with respect to the spherical metric. 
 \qed \medskip

We are interested in quasisymmetric maps $g\: W\ra W'$ between  weak tangents $W$ of $S_p$ and weak tangents $W'$ of $S_q$. Note that $0,\infty\in W,W'$.  We call $g$ {\em normalized} if 
$g(0)=0$ and $g(\infty)=\infty$.

\begin{lemma} \label{L:inducmap}
Let $p,q\ge 3$ be odd, $z_0\in D_p$, $w_0\in D_q$ and $f\: S_p\ra S_q$ be a quasisymmetric map with $f(z_0)=w_0$. Then 
$f$ induces a normalized quasisymmetric map $g$ between the 
weak tangent $W$ of $S_p$ at $z_0$ and the weak tangent $W'$ of $S_q$ at $w_0$.
\end{lemma}

\no 
{\em Proof.} By   Proposition~\ref{P:Qcext} we can  extend $f$ to a quasiconformal homeomorphism $F\: \widehat \C\ra  \widehat \C$.  By our  discussion  earlier in this section, 
there exists a relative neighborhood $N$ of $z_0$ in $S_p$ and 
a  relative neighborhood $N'$ of $w_0$ in $S_q$ such that  
$W\setminus \{\infty\} = \bigcup_{n\in \N_0} p^n(N-z_0)$ and $W'\setminus \{\infty\} = \bigcup_{n\in \N_0} q^n (N'-w_0)$.
Moreover, by \eqref{eq:asc} we may assume $p^{-n}(N-z_0)\sub N-z_0$ and 
 $q^{-n}(N'-w_0)\sub N'-w_0$ for all $n\in \N_0$.

Pick a point $u_0\in N-z_0$, $u_0\ne 0$. 
 Then for each $n\in \N_0$ we have 
 $$z_0+p^{-n} u_0\in N\setminus 
 \{z_0\}\sub S_p, $$ and so $F(z_0+p^{-n} u_0)\ne  w_0,  \infty$. Hence we can 
choose a unique number $k(n)\in \Z$ as follows: if we define the map
$F_n\: \widehat \C\ra  \widehat \C$ by
$$ F_n(u)=   q^{k(n)} \big(F(z_0+p^{-n} u)-w_0\big)$$
for $u\in \widehat \C$, then 
$$ 1\le  |F_n(u_0)|<  q. $$
Note that $k(n)\to \infty $ as $n\to \infty$. Since $F(\infty)\not\in S_q$, and so $F(\infty)\ne w_0$, this implies that $F_n(\infty)\to \infty$ as $n\to \infty$. We also have $F_n(0)=0$. 
So the images of $0$, $\infty$, and $u_0\ne 0,\infty$ under $F_n$ have mutual spherical distance uniformly bounded from below independent of $n$. 
Moreover, each map $F_n$ is obtained from $F$ by pre- and post-composing by M\"obius transformations. Hence the sequence $(F_n)$ is uniformly  quasiconformal, and it follows that we can find a subsequence of $(F_n)$ that converges uniformly on  $\widehat \C$ to a quasiconformal map 
$F_\infty$~\cite[Theorem~5.1(3), p.~73]{LV}. 
 By passing to yet another subsequence if necessary, we can also assume that we have uniform convergence of the inverse maps 
 in the subsequence to $F_\infty^{-1}$.

In this way,  we can find  sequences $(k_n)$ and $(l_n)$ of natural numbers with $k_n\to \infty$ and $l_n\to \infty $ as $n\to \infty$ 
such that if we define 
$$\widetilde F_n(u)=  q^{k_n} \big(F(z_0+p^{-l_n} u)-w_0\big)$$
for $u\in \widehat \C$, then 
$\widetilde F_n\to F_\infty $ and $\widetilde F_n^{-1}\to 
F_\infty^{-1}$ uniformly on $\widehat \C$ as $n\to \infty$.  
Then  $F_\infty(0)=0$ and $F_\infty(\infty)=\infty$. Moreover, since 
$F_\infty$ is quasiconformal, this  map is a  quasisymmetry on $\widehat \C$~\cite[Theorem~4.9]{HK98}.

So to prove the statement of the lemma, it suffices to show that $F_\infty(W)=W'$, because then $g:=F_\infty|W$ is an induced normalized quasisymmetric map between  $W$ and $W'$ as desired.

Let $u\in W$ be arbitrary. If $u=\infty$, then $F_\infty(u)=\infty\in W'$. If $u\in W\setminus\{\infty\}$, then $u\in p^m(N-z_0)$
for some $m\in \N_0$, and so 
 $$z_0+p^{-l_n}u\in  z_0+(N-z_0)=N\sub S_p$$ for large $n$.   Since 
$z_0+p^{-l_n}u\to z_0$ as $ n\to \infty$, it follows that 
$$F (z_0+p^{-l_n}u) \sub N'$$ for large $n$, and so 
$\widetilde F_n(u) \in W'$. Since $W'$ is closed,
we have $$F_\infty(u)=\lim_{n\to \infty} \widetilde F_n(u)\in W'. $$ Hence $F_\infty(W)\sub W'$.

Note that 
$$\widetilde F_n^{-1}(w)= p^{l_n} \big(F^{-1}(w_0+q^{-k_n} w)-z_0\big)$$ for $w\in \widehat \C$. So we can apply the same argument to the inverse maps, and conclude  that 
$F_\infty^{-1}(W')\sub W$. It follows that $F_\infty(W)=W'$ as desired. \qed\medskip

The previous lemma is an instance of the more general fact  that  a quasisymmetric map between two  standard carpets induces a normalized  quasisymmetric maps between 
weak tangents. It is likely that such a map only exists if the weak 
tangents are isometric. If this were  the case, then this  would put strong restrictions on the 
original quasisymmetric map. Unfortunately, we are only able to prove one result in this direction.

\begin{proposition}\label{P:Corwtangents} Let $p\ge 3$ be odd. 
Then there is no normalized quasisymmetric map from $W_{{\pi}/2}(p)$ onto $W_{{3\pi}/2}(p) $.  
\end{proposition}

The proof will occupy the rest of the section. We fix an odd number $p\ge 3$
in the following. Weak tangents refer to $S_p$, and we write $W_{{\pi}/2}=W_{{\pi}/2}(p)$, etc.

We cannot prove that there is no normalized  quasisymmetric map between $W_{{\pi}/2}$ and  $W_{{\pi}} $ or between 
$W_{3{\pi}/2}$ and $W_{{\pi}}$. If this were true, the proofs of Theorems~\ref{T:Onethird} and \ref{T:Standard} would admit some  simplifications.

Let $G$  and $\widetilde G$ denote the  group    of normalized orientation-preserving quasisymmetric self-maps of $W_{{\pi}/2}$  and $W_{{3\pi}/2}$, respectively.  
 Then $G$ and $\widetilde G$  both  contain the map $z\mapsto \mu(z):=pz$ 
induced by multiplication by $p$, and so  it follows from Lemma~\ref{L:Wunifqcirc} and Corollary~\ref{C:Cycl} that $G$ and 
$\widetilde G$ are infinite cyclic.  Let $\phi$ be a generator of $G$.  It is actually very likely that $G$ is generated by $\mu$, and that we can take $\phi=\mu$, 
but there seems to be no easy proof for this statement. So the subsequent argument cannot rely on this which  causes some  complications. 
In any case,  we have $\mu=\phi^s$ for some $s\in \Z\setminus\{0\}$. By replacing $\phi$ by $\phi^{-1}$ if necessary, we may assume that $s>0$.  
 
By Lemma~\ref{L:Wunifqcirc} and  Proposition~\ref{P:Qcext} there exists a quasiconformal map
$\Phi\: \widehat \C\ra \widehat \C$ whose restriction to $W_{{\pi}/2}$
is equal to $\phi$. Then $\Phi(0)=0$ and $\Phi(\infty)=\infty$. 
Let 
 $$\Omega:=\{z\in\C\co \Re (z)>0 \text{ and } \Im (z)>0\}.$$  
 Then $C_0:=\partial \Omega$ is a peripheral circle of $W_{{\pi}/2}$ and we have $\Phi(\partial \Omega)=\phi(\partial \Omega)=\partial \Omega$. Since 
 $\phi$, and hence also $\Phi$, is orientation-preserving, 
these maps  fix the positive real axis and the positive imaginary axis setwise. It follows that  $\Phi(\Omega)=\Omega$.

Let $\Ga$  be the family of all open paths 
in the region  $\Omega$  that connect  the positive real and positive imaginary axes. By what we have just   seen, the path family $\Ga$ is $\Phi$-invariant.  The peripheral circles of $W_{{\pi}/2}$ that meet some  path in $\Gamma$ are precisely the peripheral circles $C\ne C_0=\partial \Omega$ (this is why we chose $\Gamma$ to consist of open paths). 
It follows from Corollary~\ref{C:Ccp} that $\phi^n(C)\ne C$ 
for all $n\in \Z\setminus\{0\}$ and all peripheral circles
$C$ of $W_{{\pi}/2}$ that  meet some  path in $\Gamma$.
So we can apply  Lemma~\ref{L:Cyclgp} and conclude that 
\begin{equation}\label{eq:modphimu}  \Mod_{W_{{\pi}/2}/\langle \mu\rangle}(\Gamma)= \Mod_{W_{{\pi}/2}/\langle \phi^s \rangle}(\Gamma)  =s\Mod_{W_{{\pi}/2}/G}(\Gamma). 
\end{equation}

\begin{lemma}\label{L:0infty} We have 
$0<\Mod_{W_{{\pi}/2}/G}(\Gamma)<\infty.$
\end{lemma}

\no {\em Proof.}  By \eqref{eq:modphimu} it is enough to show that
$$0<\Mod_{W_{{\pi}/2}/\langle \mu\rangle}(\Gamma)<\infty.$$

To establish the inequality  $\Mod_{W_{{\pi}/2}/\langle \mu\rangle}(\Gamma)<\infty$,  it 
suffices to exhibit  an admissible mass distribution of finite mass. 

If $C\ne C_0=\partial \Omega$ is a peripheral circle of $W_{{\pi}/2}$, we denote by $\theta(C)$ the angle under which $C$ is seen from the origin, i.e., $\theta(C)$ is the length of the circular arc obtained as the image of $C$ under the radial projection map $z\in \C\setminus \{0\}\mapsto {\rm pr}(z):= z/|z|$ to the unit circle. We set $\rho(C_0):=0$,  and 
$\rho(C):=\frac2\pi\theta(C)$ for all peripheral circles $C\ne C_0$. 
We claim that $\rho$ is admissible for $\Mod_{W_{{\pi}/2}/\langle \mu\rangle}(\Gamma)$. 

To see this,  first note that $\rho$ is constant on orbits of peripheral circles under the action of the group $\langle \mu\rangle$. Let $\Gamma_0$ denote the family of all paths $\gamma$  in $\Gamma$ that are not locally rectifiable or for which 
$\ga\cap W_{{\pi}/2}$ has positive length. Since $W_{{\pi}/2}$ is a set of measure zero, we have $\Mod(\Gamma_0)=0$. 

If $\gamma\in \Gamma$, then 
${\rm pr}(\gamma)=\alpha:=\{e^{\iu t}: 0<t<\pi/2\}$ and the projection map ${\rm pr}$  is a Lipschitz map  on $\ga$. So if $\ga\in \Gamma\setminus \Gamma_0$, then up to a set of measure zero, ${\rm pr}(\gamma)=\alpha$ is covered by the projections ${\rm pr}(C)$ of the peripheral circles that meet 
$\ga$. It follows that 
$$ \sum_{\ga \cap C\ne \emptyset} \rho(C)= \frac 2\pi
 \sum_{\ga \cap C\ne \emptyset} \theta(C)\ge 1$$
 for all $\gamma\in \Gamma\setminus \Gamma_0$. So $\rho$ is indeed admissible. 
 
To find a mass bound for $\rho$ note that every 
$\langle \mu\rangle$-orbit of a peripheral circle $C\ne C_0$ has a unique element contained in 
the set $F:=\overline {\mu(Q_0)\setminus Q_0}$ (recall that  
$Q_0=[0,1]\times [0,1]\sub \R^2\cong \C$ denotes the unit square).   Moreover, there exists a constant $K>0$ such that 
$$\theta(C)\le K \ell(C)$$ for all peripheral circles $C$ of $W_{{\pi}/2}$ with  $C\sub F$.
It follows that 
$$\mass_{W_{{\pi}/2}/\langle \mu\rangle}(\rho)=\frac 4{\pi^2}\sum_{C\sub F}\theta(C)^2\lesssim 
  \sum_{C\sub F}\ell(C)^2= {\rm Area}(F)= p^2-1,$$
  where ${\rm Area}(F)$ denotes the Euclidean area of $A$. 
  Hence $\rho $ is an admissible density for 
  $\Mod_{W_{{\pi}/2}/\langle \mu\rangle}(\Gamma)$ with finite mass as desired. 

To show that  $\Mod_{W_{{\pi}/2}/\langle \mu\rangle}(\Gamma)>0$, we argue by contradiction and assume that 
$\Mod_{W_{{\pi}/2}/\langle \mu\rangle}(\Gamma)=0$. For $k\in \N$ let $\mathcal{C}_k$ denote the set of all peripheral circles $C$ 
of $W_{{\pi}/2}$ with $C\sub F_k:=\overline{\mu^k(Q_0)\setminus \mu^{-k}(Q_0)}$. Then  every orbit $\mathcal O$ of a peripheral circle $C\ne C_0$ under the action of $\langle \mu\rangle$ has exactly $2k$ elements in common with 
$\mathcal{C}_k$. Hence $\#(\mathcal{O}\cap \mathcal{C}_k)\le N_k:=2k$. Moreover, since every path $\gamma\in \Gamma$ 
lies in $F_k$ for sufficiently large $k$, we have 
$\Gamma=\bigcup_k \Gamma_k$, where  $\Gamma_k$ 
denotes the family of all paths in $\Gamma$ that only meet peripheral circles in $\mathcal{C}_k$. This and the previous considerations imply that the hypotheses of Proposition~\ref{P:infmingroup} are satisfied. Hence there exists an extremal mass distribution for 
 $\Mod_{W_{{\pi}/2}/\langle \mu\rangle}(\Gamma)$. 
 
 By our assumption $\Mod_{W_{{\pi}/2}/\langle \mu\rangle}(\Gamma)=0$. This is only possible if every path in $\Gamma$ belongs to the exceptional family for $\Mod_{W_{{\pi}/2}/\langle \mu\rangle}(\Gamma)$. 
 We conclude that $\Mod(\Gamma)=0$;  but obviously   $\Mod(\Gamma)=\infty$, and we obtain a contradiction.
 \qed \medskip
   
Let $H$ be the group of homeomorphsims of $\widehat \C$ 
generated by the reflections in the real and in the imaginary axes. 
Then $H$ consists of precisely four elements.

 We may assume that the quasiconformal map $\Phi$ whose restriction to ${W_{{\pi}/2}}$ is equal to the generator 
 $\phi$ of $G$ has the property that it is equivariant under $H$ 
 in the sense that $\Phi\circ \alpha=\alpha\circ \Phi$ for all $\alpha\in H$. 
 
  Indeed, if the original map $\Phi$ does  not have this property, then we 
 restrict it to the first quadrant and extend this restriction by 
 successive reflections in real and imaginary axes to the whole 
 sphere. The new map $\Phi$ obtained in this way is  clearly  an
 orientation-preserving  homeomorphism with the desired equivariance property. It is also quasiconformal away from the real and positive imaginary axes. Since sets of finite $1$-dimensional Hausdorff measures form  removable singularities for quasiconformal maps~\cite[Theorem~3.2, p.~202]{LV}, 
 $\Phi$ will actually be a quasiconformal map on $\widehat \C$. As before, $\Phi|{W_{{\pi}/2}}=\phi$.
 
 Let 
 $$\widetilde \Omega:=\{ z\in \C: \Re (z)<0 \text { or } \Im (z)<0\}.$$
 Then $\widetilde  \Omega$ is a three-quarter plane whose closure contains 
 $W_{{3\pi}/2}$, and $C_0=\partial \Omega=\partial \widetilde \Omega$ is a peripheral circle of $W_{{3\pi}/2}$. The set  $W_{{3\pi}/2}$ consists of three 
  copies of $W_{{\pi}/2}$ that can be obtained by successive reflections 
in   the real and positive imaginary axes. By its  equivariance 
property the map  $\Phi$ restricts to a normalized orientation-preserving quasisymmetric 
self-map  $\psi:=\Phi|W_{{3\pi}/2}$ of 
$W_{{3\pi}/2}$.  

Recall that $\widetilde G$  denotes the infinite cyclic  group  consisting 
of all normalized  orien\-tation-preserving 
quasisymmetric self-maps of  $W_{{3\pi}/2}$. Then  we have $\psi\in \widetilde G$, and $\langle \psi \rangle$ is an infinite 
cyclic subgroup  of $\widetilde G$.  

Let  $\widetilde \Gamma$ be the family of all open paths in $\widetilde \Omega$ that join the  positive real and the positive imaginary axes. 

\begin{lemma}\label{L:serial} We have 
$\Mod_{W_{3{\pi}/2}/\langle \psi \rangle}(\widetilde \Gamma)\le \tfrac13 \Mod_{W_{{\pi}/2}/G}(\Gamma)$. 
\end{lemma}

\no 
{\em Proof.} Essentially, this follows from an application of a 
suitable ``serial law" to modulus with respect to a group. 

More precisely, suppose that $\rho$ is an arbitrary admissible 
 invariant 
mass distribution for $\Mod_{W_{{\pi}/2}/G}(\Gamma)$ with exceptional family $\Gamma_0$. 
We want to  use $\rho$ to define a suitable admissible mass distribution 
$\tilde \rho$ for $\Mod_{W_{3{\pi}/2}/\langle \psi \rangle}(\widetilde \Gamma)$. For the special peripheral circle 
$C_0$ of $W_{3{\pi}/2}$ we set $\tilde \rho(C_0)=0$.
If $ \widetilde C$ is any peripheral circle of $W_{3{\pi}/2}$ with $ \widetilde C\ne C_0$, then there exists a unique element $\alpha\in H$ such that 
$\alpha(\widetilde C)$ is a peripheral circle of $W_{{\pi}/2}$.
We  set $\tilde\rho(\widetilde C):=\frac 13 \rho( \alpha(\widetilde C))$. 

By the equivariance property of $\Phi$ and the fact that 
$\rho$ is constant on orbits of $G=\langle \phi\rangle$, it follows 
that $\tilde\rho$ is constant on orbits  of $\langle \psi\rangle$. 

Let $\widetilde \Gamma_0$ be family of all paths 
in $\widetilde \Gamma$ that have a subpath that  can be mapped to a path in $\Gamma_0$ by an element $\alpha\in H$. 
Since $\Mod(\Gamma_0)=0$, we have 
$\Mod(\widetilde \Gamma_0)=0$.

Let $\ga\in \widetilde \Gamma$ be arbitrary. Then $\gamma$ 
has three disjoint open subpaths (one for each quarter-plane
of $\widetilde \Omega$) that are mapped to a path in $\Gamma$ by a suitable element in $H$. Let $\gamma_i$, $i=1,2,3$, denote 
these image paths in $\Gamma$. If in addition  $\ga\not\in\widetilde \Gamma_0$, then   $\gamma_i\not\in \Gamma_0$ for $i=1,2,3$;  so
$$\sum_{\ga \cap \widetilde C\ne \emptyset } 
\widetilde \rho(\widetilde C)\ge \frac 13 \sum_{i=1}^3\sum_ 
{\ga_i \cap  C\ne \emptyset }  \rho(C)\ge 1 $$
for all $\ga\in \widetilde \Gamma\setminus \widetilde \Gamma_0$. 
Hence $\widetilde \rho$ is admissible 
for $\Mod_{W_{3{\pi}/2}/\langle \psi \rangle}(\widetilde \Gamma)$
and it follows that 
$$ \Mod_{W_{3{\pi}/2}/\langle \psi \rangle}(\widetilde \Gamma)
\le \mass_{W_{3{\pi}/2}/\langle \psi \rangle}(\tilde \rho)\le \tfrac 13 
 \mass_{W_{{\pi}/2}/ G}( \rho). $$
 Since $\rho$ was an arbitrary admissible mass distribution for 
$\Mod_{W_{{\pi}/2}/G}(\Gamma)$, the statement  follows. 
\qed\medskip

\noindent
\emph{Proof of Propsition~\ref{P:Corwtangents}.} We use the notation introduced above, and denote by $G$ and $\widetilde G$ infinite   cyclic groups  
of all normalized orientation-preserving 
quasisymmetric self-maps of  $W_{{\pi}/2}$ and $W_{{3\pi}/2}$,
respectively. As before 
let $\Ga$ and $\widetilde \Ga$ be the family of all paths 
in $\Omega$ and $\widetilde \Omega$, respectively,  that join the positive real and the positive imaginary axes.

To reach a contradiction, we assume  that there is a normalized quasisymmetric
map $f$ from $W_{{\pi}/2}$ onto $W_{{3\pi}/2}$. Precomposing $f$ by the reflection in the line $L=\{z\in \C: \Re (z)=\Im (z)\}$ if necessary, we may assume that  $f$ is  orientation-preserving.  Then $\widetilde G=f\circ G\circ f^{-1}$,  and $\tilde\phi:=f\circ\phi\circ f^{-1}$ is a  
generator for $\widetilde G$.  By Lemma~\ref{L:Wunifqcirc} and Proposition~\ref{P:Qcext}, the map $f$ extends to a quasiconformal map $F$  on $\widehat\C$.
Then $\widetilde\Ga=F(\Ga)$, and so 
Lemma~\ref{L:Invgpmod}  gives
$$
{\rm mod}_{W_{{3\pi}/2}/\widetilde G}(\widetilde\Ga)={\rm mod}_{W_{{\pi}/2}/G}(\Ga).
$$

Let  $\psi=\Phi|W_{{3\pi}/2}\in \widetilde G$ be the map considered above. Then $\psi=\tilde \phi^m$ for some $m\in \Z\setminus \{0\}$,  and it follows from Corollary~\ref{C:Ccp} and  Lemma~\ref{L:Cyclgp} (see the argument that we used to establish \eqref{eq:modphimu}) that 
$$  {\rm mod}_{W_{{3\pi}/2}/\langle  \psi\rangle }(\widetilde\Ga)=
 |m|\, {\rm mod}_{W_{{3\pi}/2}/\widetilde G}(\widetilde\Ga). $$
 
 Hence by Lemma~\ref{L:serial} we have 
\begin{eqnarray*}  {\rm mod}_{W_{{\pi}/2}/G}(\Ga)&=&
{\rm mod}_{W_{{3\pi}/2}/\widetilde G}(\widetilde \Ga)\\
&=& \tfrac 1{|m|}{\rm mod}_{W_{{3\pi}/2}/\langle\psi\rangle} (\widetilde \Ga)\\
&\le &\tfrac 1{3|m|} {\rm mod}_{W_{{\pi}/2}/G}(\Ga).
\end{eqnarray*}  
This is only possible if ${\rm mod}_{W_{{\pi}/2}/G}(\Ga)=0$ or 
${\rm mod}_{W_{{\pi}/2}/G}(\Ga)=\infty$;  this  contradicts Lemma~\ref{L:0infty} and  the statement follows. 
 \qed \medskip

\section{Proof of Theorems~\ref{T:Onethird}--\ref{T:Standard}}\label{S:PT1}

\no 
We fix  an odd integer  $p\ge 3$. 
As before we assume that the standard Sierpi\'nski carpet $S_p$ is obtained by subdividing the unit square $Q_0=[0,1]\times [0,1]$ in the first quadrant of $\C\cong \R^2$. In this section it is convenient to  mostly use  real notation; so $(x_0,y_0)$ is the point in $\R^2$ with  $x$-coordinate $x_0$ and $y$-coordinate $y_0$. As before we use $0$ to denote the origin in $\R^2$. 

   We equip $S_p$    with the restriction of the Euclidean metric.   The carpet $S_p$ has four lines of symmetries; one of them is the diagonal $D:=\{(x,y)\in \R^2:x=y\}$ 
and another the vertical line $V:=\{(x,y):x=1/2\}$. We denote 
the reflections in $D$ and $V$ by $R_D$ and $R_V$, respectively. The maps $R_D$ and $R_V$ generate the 
group of Euclidean isometries of $S_p$, which consists of eight elements.

If  $f$ is  a quasisymmetric self-map of $S_p$, then   
by  Corollary~\ref{C:Group}  the outer square $O$ and the middle square $M$ of $S_p$  are preserved as a pair; so $ \{f(O), f(M)\}=\{O,M\}$.  We will now show that $f(O)=M$ is actually impossible.

\begin{lemma}\label{L:f(O)=O} Let $f$ be a quasisymmetric self-map of $S_p$, $p\ge 3$ odd. Then $f(O)=O$ and $f(M)=M$.  \end{lemma}

\no{\em Proof.} 
Let $f\in {\rm QS}(S_p)$ be arbitrary.  We know that $f(O)\in \{O,M\}$.  It is enough to show
 $f(O)=O$, because then necessarily $f(M)=M$. We argue by contradiction and assume that $f(O)=M$.  Then $f(M)=O$, and so $f$ interchanges $O$ and $M$. 

By Corollary~\ref{C:Gpfinite}  the group ${\rm QS}(S_p)$ of all quasisymmetric self-maps of $S_p$  is finite. 
Let $G$ be the  subgroup of ${\rm QS}(S_p)$
consisting of all quasisymmetric self-maps $g$ of $S_p$ 
with $g(O)=O$ and $g(M)=M$. Then $G$ is  also finite and  contains the isometry group of $S_p$. Moreover, if $G_0$ is the set of all maps in $G$ that are orientation-preserving, then $G_0$ is a subgroup in  $G$ 
  of index $2$; indeed, $G$ can be written as the disjoint
  union 
  \begin{equation}\label{eq:index2}
  G=G_0\cup G_0R_D
  \end{equation}
   of two right cosets of 
  $G_0$. 
  
  If $z\in S_p$ is arbitrary, we denote by 
  $$\mathcal{O}(z)=\{g(z): g\in G\}$$ 
  the orbit of $z$ under the action of $G$.
  Let  $$c=({(p-1)}/{(2p)},{(p-1)}/{(2p)})$$ be the left lower corner of the square $M$, and $w_0=f(0)\in M$. In the following we will 
  consider the orbits $\mathcal{O}(c)$ and $\mathcal{O}(w_0)$.
  Both are subsets of $M$. 
  Since $G$ contains the isometry group of $S_p$, the orbits 
$\mathcal{O}(w_0)$ and  $\mathcal{O}(c)$ have the same  
symmetries as  $S_p$.
  
  It  follows from Corollary~\ref{C:Ccp} that if $g\in G_0$ has a fixed point in $S_p$, then $g$ is the identity on $S_p$.
 So if $z\in S_p$ is arbitrary, then the map $g\in G_0\mapsto g(z)$ is injective. Since $R_D(0)=0$, it follows from \eqref{eq:index2}
 that 
 $\#\mathcal{O}(0)=\#G_0.$
 
  We also have $R_D(c)=c$, and so $\#\mathcal{O}(c)=\#G_0.$
Moreover, the map $g\in G\mapsto f\circ g \circ f^{-1} \in G$ is an automorphism of $G$. This implies that 
$$\mathcal{O}(w_0)=\{(f\circ g \circ f^{-1})(w_0): g\in G\}=\{(f\circ g) (0): g\in G\} =f(\mathcal{O}(0)). $$
Hence 
$$ \#\mathcal{O}(w_0)=\#\mathcal{O}(0)=\#G_0= \#\mathcal{O}(c), $$
and so the orbits $\mathcal{O}(w_0)$ and  $\mathcal{O}(c)$
have the same number of elements.

We will now show that this is impossible. 
First note that $c\not\in \mathcal{O}(w_0)$, and so $\mathcal{O}(c)\cap  \mathcal{O}(w_0)=\emptyset$. Indeed, suppose on the contrary that $c\in \mathcal{O}(w_0)$. Then there 
exists $g\in G$ with  $c=g(w_0)=(g\circ f)(0)$. Then 
$h:=g\circ f$ is a quasisymmetric self-map of $S_p$ with 
$h(0)=c$. By Lemma~\ref{L:inducmap} the map $h$ induces a normalized quasisymmetric map from  
$W_{\pi/2}$, the weak tangent of $S_p$ at $0$,  onto $W_{3\pi/2}$,
the weak tangent of $S_p$ at $c$.  This is impossible by Proposition~\ref{P:Corwtangents}. 

By symmetry $\mathcal{O}(c)$ contains all corners of $M$, while $\mathcal{O}(w_0)$ contains none of the 
corners of $M$ by what we have just seen. 

Let 
$$m'=(1/2,(p-1)/(2p))$$
be the midpoint of the bottom  side of $M$.  We want to show that 
$m'$ belongs to neither $\mathcal{O}(w_0)$ nor  $\mathcal{O}(c)$. Indeed, suppose that $m'\in \mathcal{O}(w_0)$. Similarly as above,  we  can then
find  a 
 quasisymmetric self-map $h$ of $S_p$ with 
$h(0)=m'$. By precomposing $h$ by $R_D$ if necessary, we may 
assume that $h$ is orientation-preserving. Since the weak tangent 
of $S_p$ at $m'$ is isometric to $W_\pi$, we get an induced normalized quasisymmetric map $h_1\: W_{\pi/2}\ra W_\pi$.

We necessarily have $h(O)=M$ and $h(M)=O$. Consider the map 
$R_V\circ h \circ R_D$. Then $h$ and $R_V\circ h \circ R_D$
are orientation-preserving quasisymmetric self-maps of $S_p$ that act in the same way on  the origin, and on the peripheral circles $O$ and $M$. 
By Corollary~\ref{C:Ccp}  
it follows that $h=R_V\circ h \circ R_D$.  This shows that $h$ maps 
the set $S_p\cap D$ onto $S_p\cap V$.
Since we know that $h(c)\in h(M)=O$, this only leaves two possibilities for the point $h(c)$, namely $(1/2, 0)$ and 
$(1/2, 1)$. Since at both points the weak tangents of
$S_p$ are isometric to $W_\pi$, we get an induced normalized quasisymmetric map $h_2\: W_{3\pi/2}\ra W_\pi$.
Then $h_2^{-1}\circ h_1$ is a normalized quasisymmetric map from  $W_{\pi/2}$ onto $W_{3\pi/2}$. We get a contradiction to 
Proposition~\ref{P:Corwtangents}, showing that $m'\not\in 
\mathcal{O}(w_0)$. 

The proof that $m'\not\in 
\mathcal{O}(c)$ runs along similar lines. Again we argue by contradiction and assume $m'\in 
\mathcal{O}(c)$. Then we can find an orientation-preserving quasisymmetric 
self-map $h$ of $S_p$ with $h(c)=m'$. This gives a normalized  quasisymmetric map $h_1\: W_{3\pi/2}\ra W_{\pi}$.
We must have $h(M)=M$ and $h(O)=O$. Then $h$ and $R_H\circ  h \circ R_D$ are orientation-preserving quasisymmetric self-maps of $S_p$ that act in the same way on $c$, and on the peripheral circles $O$ and $M$. Therefore, 
$h=R_V\circ h \circ R_D$, and so $h$ maps  
the set $S_p\cap D$  onto $S_p\cap V$.
This only leaves the possibilities $(1/2, 0)$ or $(1/2, 1)$ 
for the point $h(0)$. In any case, 
we get an induced normalized quasisymmetric map $h_2\: W_{\pi/2}\ra W_\pi$, and by considering  $h_1^{-1}\circ h_2$,  again   a contradiction to 
Proposition~\ref{P:Corwtangents}. Hence  $m'\not\in 
\mathcal{O}(c)$.

To summarize, we know that the sets $\mathcal{O}(c)$ and 
$\mathcal{O}(w_0)$ are disjoint subsets of $M$ with   the same symmetries as $S_p$,  and none of these sets contains $m'$. This implies that each side of $M$ contains an even number of points in  $\mathcal{O}(c)$ and 
$\mathcal{O}(w_0)$, since we have reflection symmetry about the midpoint of each side.

Since  no corner of $M$ is in $\mathcal{O}(w_0)$  and each  side of $M$ contains the same even number of points in $\mathcal{O}(w_0)$, it follows that 
$$\#\mathcal{O}(w_0)=8k$$
for some $k\in \N_0$.  So $\#\mathcal{O}(w_0)$ is divisible by 
$8$. On the other  hand, $\mathcal{O}(c)$ contains the corners 
of $M$. Since each corner belongs to two sides, we have  $$\#\mathcal{O}(c)=8l-4$$
for some $l\in \N$, and so $\#\mathcal{O}(c)$ is not divisible by $8$. Since we know that  $\#\mathcal{O}(w_0)=\#\mathcal{O}(c)$,
this is a contradiction. So  $f(O)=M$ is impossible, and we must have $f(O)=O$. 
\qed\medskip

One can make the logic of the previous proof a little more transparent, if one follows a slightly different (albeit longer) route. Namely, if a map $f\in {\rm QS}(S_p)$
with $f(O)=M$ exists, then, by using a counting argument as above, one can find such a map $f$ that sends the origin to one of the natural candidates adapted to the symmetries of $S_p$, namely   to  a corner of $M$ or to  the midpoint of one of the sides of $M$. Arguing as in the previous proof based on Proposition~\ref{P:Corwtangents}, one can rule out these possibilities, and again reaches a contradiction.

As before let  $O$ be the outer and $M$ the middle square of $S_p$. We denote the orbit of a point $z\in S_p$ by the
group ${\rm QS}(S_p)$ of quasisymmetric self-maps of $S_p$
by $\mathcal{O}(z)$. Now that we know that every map 
$f\in {\rm QS}(S_p)$ preserves $O$ and $M$ setwise, 
the group $G$ introduced in the previous proof is actually 
equal to ${\rm QS}(S_p)$. Define $m=(0,1/2)$. Then $m$ is the midpoint of the bottom side of $O$. 

\begin{lemma}\label{L:8div} Let $z\in O$ be arbitrary. 
If $\mathcal{O}(z)\ne \mathcal{O}(0), \mathcal{O}(m)$, then 
$\#\mathcal{O}(z)$ is divisible by $8$. 
Moreover,   $\mathcal{O}(0)$ and $\mathcal{O}(m)$
are divisible by $4$, but not by $8$, and $\mathcal{O}(0)\cap 
\mathcal{O}(m)=\emptyset$. 
 \end{lemma}
 
 \no {\em Proof.} If $z\in O$, then the  orbit $\mathcal{O}(z)\sub O$ has  the same symmetries as $S_p$; so each side of $O$ contains the same number
 of points in  $\mathcal{O}(z)$, and $\#\mathcal{O}(z)$ must be divisible by $4$.  If $\mathcal{O}(z)\ne \mathcal{O}(0), \mathcal{O}(m)$, then $\mathcal{O}(z)$ does not contain any corners  of $O$, nor any midpoint of a side of $O$. Hence 
 each side of $O$ contains an even number $2k$, $k\in \N$, of points in $\mathcal{O}(z)$ and $\#\mathcal{O}(z)=8k$.  In this case, $\#\mathcal{O}(z)$ is divisible by $8$.
 
 We want to show that $\mathcal{O}(0)\cap 
\mathcal{O}(m)=\emptyset$. We argue by contradiction, and assume that  $\mathcal{O}(0)\cap 
\mathcal{O}(m)\ne \emptyset$. Then we can find a map $f\in {\rm QS}(S_p)$ with $f(0)=m$, and pre-composing $f$ with $R_D$ if necessary, we may assume that $f$ is orientation-preserving. Similarly as in the proof of 
Lemma~\ref{L:f(O)=O} this leads to a contradiction; 
namely, we first get an induced normalized quasisymmetry $f_1\: 
W_{\pi/2}\ra W_{\pi}$. 
Moreover, from  
Corollary~\ref{C:Ccp} we conclude that $f=R_V\circ f\circ R_D$, and so $f$ maps $S_p\cap D$ onto $S_p\cap V$; hence the lower 
left corner $c$ of $M$ must be mapped to the intersection $M\cap V$ which consists of two points where the weak tangent is isometric to $W_{\pi}$. This gives an induced normalized quasisymmetry  $f_2\:W_{3\pi/2}\ra W_{\pi}$. Considering $f_2^{-1}\circ f_1$ we get a contradiction from Proposition~\ref{P:Corwtangents}. 

It follows that  $\mathcal{O}(0)$ contains the corners of $O$, but not any midpoint of a side. So the number of points in 
$\mathcal{O}(0)$ on each side is an even number $2r$, $r\in \N$, and we conclude $\#\mathcal{O}(0)=8r-4$.

Finally, the number of points in $\mathcal{O}(m)$ on each side is an odd number $2l-1$, $l\in \N$, since $\mathcal{O}(m)$
contains  the midpoint of the  side. Since none of the corners of 
$O$ belongs to $\mathcal{O}(m)$, we have 
$\#\mathcal{O}(0)=8l-4$. Hence neither   
$\#\mathcal{O}(0)$ nor $\#\mathcal{O}(m)$ is divisible by $8$. 
\qed \medskip
  
 Exactly the same statement with essentially the same proof  is true for orbits of points in $M$, if in Lemma~\ref{L:8div} we replace $0$ by a corner of $M$, and  $m$ by a midpoint of a side of  $M$.

\medskip\no
{\em Proof of Theorem~\ref{T:Onethird}.} 
Let $f$ be a quasisymmetric self-map of $S_3$. We want to show that $f$ is a Euclidean  isometry of $S_3$. To see this,  we may assume $f$ is orientation-preserving, for otherwise we can compose this map   with the  reflection $R_V$ which lies in the isometry group of $S_3$.
By   Lemma~\ref{L:f(O)=O} we know that if $O$ is  the outer and $M$  the middle square of $S_3$,
then $f(O)=O$ and $f(M)=M$.   

There are eight peripheral circles of $S_3$ that are squares of sidelength  $1/9$. We call them {\em second generation squares} as they are the boundaries of the solid squares that were removed in the second step of the construction (in the first step, the square bounded by $M$ was removed from $Q_0$).   Four second generation squares, the  {\em corner squares},  have distance $1/9$ to precisely 
two sides of the unit square $Q_0$; the four other ones, the {\em side squares},  have distance $1/9$ to exactly one 
side of $Q_0$.  

Before continuing, we give a general outline of the 
ensuing argument. 
We will show that $f$ must map some  second generation 
square to another one. This will lead to various combinatorial possibilities. We will analyze them in detail. In some cases we can invoke the Three-Circle Theorem, Corollary 
\ref{C:TCT}, to identify $f$ with an isometry as desired. 
In the other cases, a  map with the given mapping behavior 
on  a second generation square will not exist.
The strategy for ruling out the existence of such ``ghost maps" is this: using symmetries and again the 
Three-Circle Theorem, we will be able to restrict  the possibilities 
for the image of the origin under $f$. Once we know that $f(0)=p\in S_3$, the map $f$ will induce a quasisymmetry
of the weak tangent $W_{\pi/2}$ of $S_3$ at $0$ to a weak tangent of $S_3$ at $p$. As we will see, this always  leads to a normalized quasisymmetry from  $W_{\pi/2}$ onto $W_{3\pi/2}$. Invoking 
Proposition~\ref{P:Corwtangents} we will then get  a contradiction ruling out the existence of the  map.  

We now proceed 
to presenting the details.

\smallskip \no 
{\em Claim.} The map $f$ sends some  second generation square 
to another second generation square. 
\smallskip

Among the eight second generation   squares let  $C_0$ be one for which 
$
{\rm mod}_{S_3}(\Ga(C_0,O;S_3))
$
is largest, and define 
 $C_1=f(C_0)$. 
Then   $C_1$ is a peripheral circle of $S_3$ and hence a square. 
As in the proof of Corollary~\ref{C:Group}, Lemma~\ref{L:Invmod} implies that
$$
{\rm mod}_{S_3}(\Ga(C_0,O;S_3))={\rm mod}_{S_3}(\Ga(C_1,O;S_3)). 
$$
For establishing the claim it suffices to show that $C_1$ has sidelength $1/9$ and is hence a second generation square.
Since $f(O)=O$ and $f(M)=M$, we have 
$C_1=f(C_0)\not\in \{O,M\}$ and so the sidelength of $C_1$ is at most   $1/9$. The 
 monotonicity of carpet  modulus and the self-similarity of $S_3$ imply that  this side length cannot be strictly smaller than $1/9$.  Indeed, suppose that this is the case.
Then there exists a unique carpet $S\sub S_3$  that can be mapped to $S_3$ by a Euclidean similarity so that 
  $C_1$ corresponds to  second generation square of $S_3$. Let $o$ denote the outer peripheral circle of $S$ corresponding to $O$. By our assumption that $C_1$ is not a second generation square, we have $S\ne S_3$,  and so $S$ is a proper subset of $S_3$.

 By definition of $C_0$ and scale invariance of carpet modulus we have 
$$
{\rm mod}_S(\Ga(C_1,o;S))\leq{\rm mod}_{S_3}(\Ga(C_0,O;S_3)).
$$ 
 On the other hand, an argument  as in the proof of 
 Lemma~\ref{L:Pair} gives
$$
{\rm mod}_{S_3}(\Ga(C_1,O;S_3))<{\rm mod}_S(\Ga(C_1,o;S)).
$$ 
The previous three  modulus relations  combined lead to a contradiction, and the claim follows. 

\smallskip
Having established that the the image $C_1=f(C_0)$ of some second generation square $C_0$ is also a second generation square, we now distinguish several cases depending on the type of the squares $C_0$ and $C_1$, i.e., whether they are corner or sides squares.  These cases will  exhaust all possibilities. 

\smallskip
\no {\em Case 1:}  $C_0$ and $C_1$  are  corner squares.\smallskip 

Then   there exists an isometry $T$ of $S_3$ given by a rotation 
that maps  $C_0$ to $C_1$. 
Then $f$ and $T$ acts in the same way  on three peripheral circles  $O$, $M$, $C_0$ of $S_3$. Since $f$ and $T$ are orientation-preserving, it follows from  Corollary~\ref{C:TCT} that 
$f=T$. Hence $f$ is an isometry of $S_3$.

\smallskip 
\no {\em Case 2.} $C_0$ is a corner square, and $C_1$ is a side square. 

\smallskip 
By pre- and postcomposing $f$ by suitable rotations, we may  assume that  $C_0$ is the corner square that has distance $1/9$ to both the $x$- and $y$-axes, and that $C_1$  is the side square that has distance $1/9$ to the $x$-axis.

Then $f$ and $R_V\circ f \circ R_D$ are orientation-preserving quasisymmetric self-maps of  $S_3$  that act in the same way on $O$, $M$, $C_0$.  Again Corollary~\ref{C:TCT} allows us to conclude that $f=R_V\circ f \circ R_D$, and so $f(D\cap S_3)=V\cap S_3$.  This only leaves 
two possibilities for the image of the origin under $f$, namely the points 
$(\frac12,0)$ or  $(\frac12,1)$; since the weak tangents of $S_3$ at these points   are  isometric to $W_{\pi}$, we get an induced normalized quasisymmetric map $f_1\colon W_{{\pi}/2}\to W_{\pi}$.  

Moreover, the lower left corner  $c=(\frac13,\frac13)$ is mapped by $f$ to either $(\frac12,\frac13)$ or $(\frac12,\frac23)$. The weak tangent of $S_3$ at $(\frac13,\frac13)$ is equal to $W_{3\pi/2}$, and the weak tangents  at both points 
$(\frac12,\frac13)$ or $(\frac12,\frac23)$ are isometric to 
$W_\pi$. 
As above, the map $f$ induces a normalized quasisymmetric  map $f_2\colon W_{{3\pi}/2}\to W_{\pi}$.
Then the  map $f_2^{-1}\circ f_1$ is a normalized quasisymmetric map from
$W_{{\pi}/2}$ onto $W_{{3\pi}/2}$. This  contradicts 
Proposition~\ref{P:Corwtangents}, and a map $f$ as in this case does not exist. 

\smallskip 
\no {\em Case 3.} $C_0$ is a side square, and $C_1$ is a 
corner  square. 
\smallskip 

Then we consider $f^{-1}$, and reduce to Case 2. This shows that a map $f$ as in this case does not exist.

\smallskip 
\no {\em Case 4.} $C_0$ and $C_1$ are side squares.\smallskip 

Then there exists a rotation $T$ of $S_3$ that maps $O$, $M$, 
$C_0$ in the same way as $f$, and as in Case 1, we conclude that 
$f=T$. Hence $f$ is an isometry of $S_3$.

The Cases~1.--4.\    exhaust all possibilities, and we have shown that in each case the map $f$ is an isometry or does not exist. Theorem~\ref{T:Onethird} follows.  
\qed \medskip



\begin{rems}\label{R:Onefifth}
\rm{By using a  similar technique and a slightly more refined case  analysis, one can show the same  
 rigidity result for the standard Sierpi\'nski carpet $S_5$; 
namely, every quasisymmetric self-map of $S_5$ is a Euclidean  isometry. 
For larger  $p$ the number of second generation squares of $S_p$ increases  and the case analysis seems to meet insurmountable obstacles.}

\rm{A more natural approach is  to first prove rigidity statements for weak tangents of $S_p$. In view of   Theorem~\ref{P:Cor} or Theorem~\ref{P:ACyl} one may speculate whether 
a normalized quasisymmetry between 
two weak tangents of a carpet $S_p$ only exists if the weak tangents are similar, i.e., one is the image of the other by a Euclidean similarity. If this is the case,  
then by considering  $W_{\pi/2}$, one can conclude that under any quasisymmetry $f$ of $S_p$, the origin must be mapped to  a
 corner of the unit square,  
 and it would easily  follow from Corollary~\ref{C:Ccp} that 
 $f$ is an isometry of $S_p$.
 
 Unfortunately,  we cannot even rule out   the  existence of a normalized quasisymmetric map between the  weak 
 tangents   $W_{\pi/2}$ and 
 $W_{\pi}$ of $S_p$. This caused some complications in the previous proof that we were able to overcome  by ad hoc arguments.} 
\end{rems}

\medskip \no
{\em Proof of Theorem~\ref{T:Onepth}.} Let $G={\rm QS}(S_p)$ be the group of all quasisymmetric self-maps of $S_p$, and $G_0$ be the subgroup of all orientation-preserving maps in $G$. 
Then $G_0$ is a subgroup in $G$ of index $2$, and $G_0$ is finite cyclic as follows from  
Lemma~\ref{L:f(O)=O} and Corollary~\ref{C:CC}. 

Consider the orbit $\mathcal{O}(0)$  of the origin under $G$. Since $R_D(0)=0$, this set is equal to the orbit of $0$ under $G_0$. Since each element 
in $G$ preserves the outer square $O$, we know that $\mathcal{O}(0)$
consists of points on  $O$. Moreover,
$\mathcal{O}(0)$ is symmetric with respect to all symmetries of $S_p$. We equip $O$ with positive orientation, so that 
$S_p$ lies on the left if we run through $O$ with this orientation.
Let $z_0=0, z_1, \dots, z_{n-1}, z_{n}=z_0$, where $n\in \N$,   be the points in $\mathcal{O}(0)$
in cyclic  order on the oriented curve $O$, and let $\alpha_i$
for $i=0, \dots, n-1$ be the corresponding subarcs of $O$ with endpoints 
$z_{i-1}$ and $z_{i}$. 
 
 There exists an element $r\in G_0$ with $r(z_0)=z_1$. 
Then $r(\alpha_0)$ is a subarc of $O$  that has the initial point  $z_1$, is positively oriented on $O$, since $r$ is orientation-preserving, and has 
its endpoint in $\mathcal{O}(0)$. Moreover, $r(\alpha_0)$ does not contain any point from $\mathcal{O}(0)$ in its interior, because this is true for $\alpha_0$. Hence the endpoint of $r(\alpha_0)$ must be $z_1$ and so $r(\alpha_0)=\alpha_1$. Repeating this argument successively for the arcs $\alpha_1,
\dots, \alpha_{n-1}$, we conclude that $r(\alpha_i)=\alpha_{i+1}$ for all $i=0, \dots, n-1$, where $\alpha_{n}=\alpha_0$. 
In particular, $r(z_i)=z_{i+1}$ and so $z_i=r^i(0)$ for 
$i=0, \dots, n$. 

This implies that $r$ generates $G_0$; indeed, if $g\in G_0$ is 
arbitrary, then by what we have just seen, there exists 
$i\in \{0,\dots, n-1\}$ such that $g(0)=z_i=r^i(0)$. 
Then $g^{-1}\circ r^{i}$ is an orientation-preserving element 
in $G$ that fixes the origin, and the peripheral   circles $O$ and $M$. Hence $g^{-1}\circ r^{i}=e$, where $e={\rm id}_{S_p}$, and so $g=r^i$. Since $r^n(0)=r^n(z_0)=z_{n}=z_0=0$, the same argument shows that $r^n=e$. Moreover, since the points 
$z_i=r^i(0)$ for $i=0, \dots, n-1$ are all distinct, $n$ is the order of
$g$. 

 Let $s=R_D$ be the reflection 
in $D$. Then $s\in G$ is orientation-reserving, and, since $G_0$ has index $2$ in $G$, it follows that $s$ and $r$ generate 
$G$. Since the orbit $\mathcal{O}(0)$ is invariant under $s$, the arc $s(\alpha_0)$ has its endpoints in $\mathcal{O}(0)$. 
There are no points from the orbit in its interior, one of the endpoints is $z_0=0$, and $s(\alpha_0)$ is traversed in negative orientation if we traverse $\alpha_0$ positively. Hence 
$s(\alpha_0)=\alpha_{n-1}$, and so $s(z_1)=z_{n-1}$.
It follows that $(s\circ r)^2\in G$ is orientation-preserving
and 
\begin{eqnarray*} (s\circ r)^2(0)&=&(s\circ r\circ s) (r(0))\,=\,(s\circ r\circ s)(z_1)\\ &=&
(s\circ r)(z_{n-1})\,=\, s(z_0)=s(0)=0.
\end{eqnarray*}
Similarly as before,  we conclude that $(s\circ r)^2=e$.
So we have the relations $s^2=r^n=(s\circ r)^2=e$ for the generators $s$ and $r$ of $G$. Moreover, 
the element $r$ has order $n$. The elements   $s$ and $s\circ r$ are orientation-reversing, and so they  have order $2$.   This implies  that $G$ is a finite dihedral group. 
\qed\medskip

\begin{rem}\label{rem:oricyc} Let $G_0$ be  the group of orientation-preserving maps in $G={\rm QS}(S_p)$.  As we have seen in the preceding proof,  $G_0$ is a cyclic subgroup of $G$ with index $2$. Moreover, the order $n$ of $G_0$ is equal to the cardinality of the orbit of $\mathcal{O}(0)$ of $0$ under 
$G$. So by Lemma~\ref{L:f(O)=O} the order $n$ of $G_0$ is divisible by $4$, but not by $8$. Of course, if our conjecture is true that every element in $G$ is an isometry, then $G_0$ consists 
of four rotations and $n=4$. 
\end{rem}

\no {\em Proof of Theorem~\ref{T:Standard}.} Let $p,q\ge 3$ 
be odd integers, and suppose that there exists a quasisymmetric map $f\: S_p\ra S_q$. We want to show that $p=q$. 

In the following  we use the subscript $p$ in our notation  if we refer to objects related to $S_p$ and  
$q$ if we refer to $S_q$. So $M_p$ denotes the middle square of $S_p$, etc. 

Let $G_p$ and $G_q$ be the groups of quasisymmetric self-maps of $S_p$ and $S_q$, respectively. Note that by Corollary~\ref{C:Gpfinite} the groups $G_p$ and $G_q$ are finite, and, $f$ being quasisymmetric,  conjugates $G_p$ and $G_q$.
This implies that if  $\mathcal{O}_p$ denotes the orbit of $0$ under $G_p$, then $\mathcal{O}_q:=f(\mathcal{O}_p)$
is the orbit of $f(0)$ under $G_q$.

It follows from Lemma~\ref{L:Invmod} and Lemma~\ref{L:Pair}  that $f$ maps the pair $\{O_p,M_p\}$ consisting of the outer square and  middle square of $S_q$ 
to the corresponding pair $\{O_q,M_q\}$. Hence $f(O_p)=O_q$ 
or $f(O_p)=M_q$, and in particular $f(0)\in O_q$ or $f(0)\in M_q$.

By Lemma~\ref{L:8div} the number $\#\mathcal{O}_p$ is not divisible by $8$. Since $\mathcal{O}_q$ has the same cardinality 
as  $\mathcal{O}_p$, the number $\#\mathcal{O}_q$ is not divisible by $8$ either. Applying Lemma~\ref{L:8div} and 
the remark after this lemma, we conclude that the orbit 
$\mathcal{O}_q$ of $f(0)\in M_q\cup O_q$ under $G_q$ must be equal to
the orbit of a corner of $O_q$ or $M_q$, or the orbit of a midpoint 
of a side of $O_q$ or $M_q$. 

Let $$c_q= ((q-1)/(2q), (q-1)/(2q)$$ be the lower left corner 
of $M_q$, and $m=(1/2,0)$ and 
$$m'_q=(1/2, {(q-1)}/{(2q)})$$ be the midpoint of the bottom side of  $O_q$ and $M_q$, respectively. By what we have seen, $f(0)$ must belong to an orbit of  one of the four points $0,c_q,m,m_q'$ under $G_q$. By composing $f$ with a suitable element in $G_q$, we may actually assume that $f(0)\in \{0,c_q,m,m_q'\}$. 
By pre-composing $f$ with $R_D$ if necessary, we may in addition assume that $f$ is orientation-preserving. 

We are  led to four cases that we now analyze. 

\smallskip
\no {\em Case 1.} $f(0)=0$. 
\smallskip 

Then $f(O_p)=O_q$ and $f(M_p)=M_q$. The map $f^{-1}\circ 
R_D\circ f \circ  R_D$ is an orientation-preserving quasisymmetry in $G_p$,
fixes the point $0$, and the peripheral circles $O_p$ and $M_p$
setwise.
Hence this map is equal to the identity on $S_p$ which implies 
$f\circ R_D=R_D\circ f$. 
From this in turn  we conclude that $f$ fixes the point $(1,1)$.

Let  $D'$ be the line $\{(x,y)\in \R^2: x+y=1\}$ and denote the reflection in $D'$ by $R_{D'}$. 
Then the  map  $f^{-1}\circ R_{D'}\circ f \circ  R_{D'}$
is an orientation-preserving quasisymmetry in $G_p$,
fixes the point $0$, and the peripheral circles $O_p$ and $M_p$
setwise. Hence this map is the identity on   $S_p$ and so 
$f\circ R_{D'}=R_{D'}\circ f$. It follows that $f$ fixes the points 
$(0,1)$ and $(1,0)$ or interchanges them. Since $f$ is orientation-preserving, and fixes $0$ and $(1,1)$,  this map must   fix $(0,1)$ and $(1,0)$. 
So $f$ fixes all corners of the unit square. 

By Theorem~\ref{P:Cor} the map $f$ must be the identity and the carpets $S_p$ and $S_q$ the same. Hence $p=q$.

\smallskip
\no {\em Case 2.} $f(0)=m$. 
\smallskip 

Then we get an induced normalized quasisymmetry $f_1\: 
W_{\pi/2}(p)\ra W_{\pi}(q)$.  
Moreover, by an argument as in Case 1,  we have
$f\circ R_D=R_V\circ f$. This implies that $f(S_p\cap D)= S_q\cap V$. Since we also have $f(M_p)=M_q$, 
 the lower left corner $c_p$ of $M_p$ must be mapped 
to  the midpoint  of the bottom or the top side of 
$M_q$. At these  points $S_q$ has a unique weak tangent isometric to $W_{\pi}(q)$. Hence we get an induced normalized 
quasisymmetry $f_2\: 
W_{3\pi/2}(p)\ra W_{\pi}(q)$.  Considering $f_2^{-1}\circ f_1$ we get a contradiction to  Proposition~\ref{P:Corwtangents}. 
So this case is impossible. 

\smallskip
\no {\em Case 3.} $f(0)=m'_q$. 
\smallskip 

This is very similar to Case 2. We get an induced normalized quasisymmetry $f_1\: 
W_{\pi/2}(p)\ra W_{\pi}(q)$, and have $f\circ R_D=R_V\circ f$.
Since $f(M_p)=O_q$, this limits the possible image points of $c_p$ under $f$ to the midpoints of the top or bottom side of 
$O_q$. Again we get an induced normalized 
quasisymmetry $f_2\: 
W_{3\pi/2}(p)\ra W_{\pi}(q)$, and a contradiction by Proposition~\ref{P:Corwtangents}. 

\smallskip
\no {\em Case 4.} $f(0)=c_q$. 
\smallskip 

Then we get an induced  normalized quasisymmetry $f_1\: 
W_{\pi/2}(p)\ra W_{3\pi/2}(q)$. 
We also have $f\circ R_D=R_D\circ f$, and so $f(S_p\cap D)=
S_q\cap D$. Pick a peripheral circle $C\ne O_p,M_p$ of $S_p$ that is symmetric with respect to  $D$ and let $v\in D\cap C$. 
Then $v$ is a corner of the square $C$ and so $S_p$ has a weak tangent at $v$ that is isometric to 
$W_{3\pi/2}(p)$. Moreover, $C'=f(C)$ is a peripheral circle of 
$S_q$ distinct from $O_q=f(M_p)$. It contains the point 
$v'=f(v)$ that lies on $D$. Hence $v'$ is a corner of $C'$, 
and so $S_q$ has a weak tangent at  $v'$ isometric to 
$W_{3\pi/2}(q)$. We  get an induced 
 normalized quasisymmetry $f_2\: 
W_{3\pi/2}(p)\ra W_{3\pi/2}(q)$. Considering $f^{-1}_2\circ f_1$, we again get a contradiction to Proposition~\ref{P:Corwtangents}. 
\smallskip 

In sum, only Case 1 is actually possible, and we have $p=q$ as desired.
 \qed 



\end{document}